%% file: ghost_waveguide_arxiv.tex
\documentclass[12pt]{iopart}

\usepackage{paralist}
\usepackage{graphics} %% add this and next lines if pictures should be in esp format
\usepackage{epsfig} %For pictures: screened artwork should be set up with an 85 or 100 line screen
\usepackage{graphicx}  \usepackage{epstopdf}%This is to transfer .eps figure to .pdf figure; please compile your paper using PDFLeTex or PDFTeXify.
\usepackage{mathrsfs}

\usepackage[breaklinks,colorlinks=true,linkcolor=blue,citecolor=red,backref=page]{hyperref}
\newtheorem{thm}{Theorem}[section]

\newtheorem{lemma}[thm]{Lemma}

\newtheorem{remark}[thm]{Remark}

\def\debproof{ {\bf Proof.} }
\def\finproof{\hfill {\small $\Box$} \\}

\newcommand{\eqref}[1]{(\ref {#1})}

\newcommand{\ea}{\end{eqnarray}}  
\newcommand{\ba}{\begin{eqnarray}}  
\newcommand{\ean}{\end{eqnarray*}}  
\newcommand{\ban}{\begin{eqnarray*}}  

\newcommand\iint{\int \hspace*{-0.05in} \int}

\input{mathmacros.tex}

\usepackage{iopams}  

\begin{document}

%%%%%%%%%%%%%%%%%%%%%%%%%%%%%

\title[Ghost imaging in a random waveguide]{A ghost imaging modality in a random waveguide} 

\author{Liliana Borcea$\hbox{}^{(1)}$  and  Josselin Garnier$\hbox{}^{(2)}$}

\address{$\hbox{}^{(1)}$Department of Mathematics, University of Michigan, Ann Arbor, MI 48109, USA} 
\ead{borcea@umich.edu}
 \vspace{10pt}
\address{$\hbox{}^{(2)}$Centre de Math\'ematiques Appliqu\'ees, Ecole Polytechnique, 91128 Palaiseau Cedex, France}
\ead{josselin.garnier@polytechnique.edu}
 \vspace{10pt}
\begin{indented}
\item[]
\end{indented}

%%%%%%%%%%%%%%%%%%%%%%%%%%%%%
\begin{abstract}
We study the imaging of a penetrable scatterer, aka target, in a
waveguide with randomly perturbed boundary. The target is located
between a partially coherent source which transmits the wave, and a
detector which measures the spatially integrated energy flux of the
wave. The imaging is impeded by random boundary scattering effects
that accumulate as the wave propagates. We consider a very large distance
(range) between the target and the detector, where that cumulative
scattering is so strong that it distributes the energy evenly among the components (modes)
of the wave. Conventional imaging is impossible in this equipartition
regime. Nevertheless, we show that the target can be located with a
ghost imaging modality. This forms an image using the
cross-correlation of the measured  energy flux, integrated over the aperture of the detector, with the time and space
resolved energy flux in a reference waveguide, at the search range. We
consider two reference waveguides: The waveguide with unperturbed
boundary, in which we can calculate the energy flux, and the actual
random waveguide, before the presence of the target, in which the
energy flux should be measured. We analyze the ghost imaging modality
from first principles and {show that it can be efficient in a 
random waveguide geometry in which there is both strong modal dispersion 
and mode coupling induced by scattering, provided that the standard ghost imaging
function is modified and integrated over a suitable time offset window in order to 
compensate for dispersion and diffusion.} The analysis 
quantifies the resolution of the image in terms of 
the source and detector aperture, the range offset between the source
and the target, and the duration of the measurements.

\end{abstract}

\vspace{2pc}

\maketitle

%%%%%%%%%%%%%%%%%%%%%%%%%%%%%%%%%%%%%%%%%%%%%%%%%%%
\section{Introduction}
\begin{figure}[t]
\begin{center}
\vspace{0in} \includegraphics[width=0.8\textwidth]{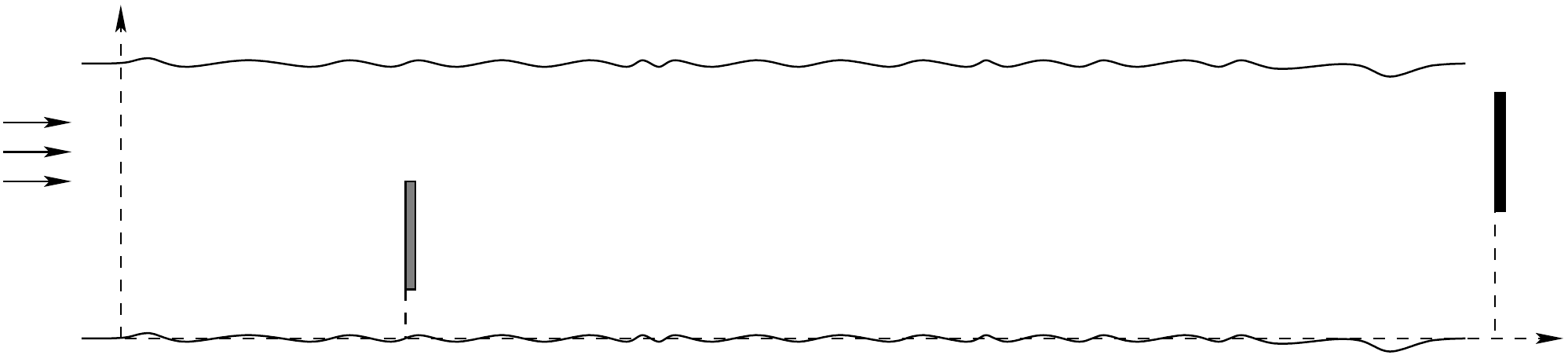}
\end{center}
\setlength{\unitlength}{3947sp}%
\begingroup\makeatletter\ifx\SetFigFont\undefined%
\gdef\SetFigFont#1#2#3#4#5{%
  \reset@font\fontsize{#1}{#2pt}%
  \fontfamily{#3}\fontseries{#4}\fontshape{#5}%
  \selectfont}%
\fi\endgroup%
\begin{picture}(9770,1095)(1568,-1569)
\put(7600,-470){\makebox(0,0)[lb]{\smash{{\SetFigFont{7}{8.4}{\familydefault}{\mddefault}{\updefault}{\color[rgb]{0,0,0}{\normalsize $z = L + \cL$}}%
}}}}
\put(3700,-470){\makebox(0,0)[lb]{\smash{{\SetFigFont{7}{8.4}{\familydefault}{\mddefault}{\updefault}{\color[rgb]{0,0,0}{\normalsize $z = L$}}%
}}}}
\put(2500,-470){\makebox(0,0)[lb]{\smash{{\SetFigFont{7}{8.4}{\familydefault}{\mddefault}{\updefault}{\color[rgb]{0,0,0}{\normalsize $(0,0)$}}%
}}}}
\put(2600,870){\makebox(0,0)[lb]{\smash{{\SetFigFont{7}{8.4}{\familydefault}{\mddefault}{\updefault}{\color[rgb]{0,0,0}{\normalsize $x$}}%
}}}}
\put(5400,900){\makebox(0,0)[lb]{\smash{{\SetFigFont{7}{8.4}{\familydefault}{\mddefault}{\updefault}{\color[rgb]{0,0,0}{\normalsize $x = \cD^+(z)$}}%
}}}}
\put(2000,630){\makebox(0,0)[lb]{\smash{{\SetFigFont{7}{8.4}{\familydefault}{\mddefault}{\updefault}{\color[rgb]{0,0,0}{\normalsize source}}%
}}}}
\put(3950,150){\makebox(0,0)[lb]{\smash{{\SetFigFont{7}{8.4}{\familydefault}{\mddefault}{\updefault}{\color[rgb]{0,0,0}{\normalsize penetrable scatterer}}%
}}}}
\put(5400,-150){\makebox(0,0)[lb]{\smash{{\SetFigFont{7}{8.4}{\familydefault}{\mddefault}{\updefault}{\color[rgb]{0,0,0}{\normalsize $x = \cD^-(z)$}}%
}}}}
\put(8100,400){\makebox(0,0)[lb]{\smash{{\SetFigFont{7}{8.4}{\familydefault}{\mddefault}{\updefault}{\color[rgb]{0,0,0}{\normalsize detector}}%
}}}}
\end{picture}%
\vspace{-0.85in}
\caption{Illustration of the imaging setup in a waveguide with
  randomly perturbed boundary. The penetrable scatterer at range $z =
  L$ is illuminated by a partially coherent source. The detector at
  range $z = L + \mathcal{L}$ measures the spatially integrated
  energy flux of the waves, as a function of time.}
\label{fig:setup}
\end{figure}
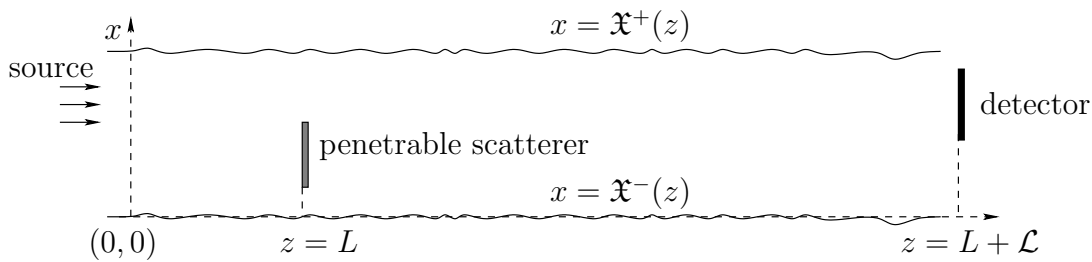

%%%%%%%%%%%%%%%%%%%%%%%%%%%%%%%%%%%%%%%%%%%%%%%%%%%
We study a ghost imaging modality in a waveguide with random
boundary. For simplicity, we consider sound waves in a two-dimensional
waveguide with straight axis and randomly perturbed sound soft
boundary, as illustrated in Figure \ref{fig:setup}. This makes it possible to
use the wave propagation theory developed in
\cite{alonso2011wave}. The results should extend qualitatively to
three dimensions, to other reflecting boundary conditions such as
sound hard \cite{alonso2011wave}, to open (radiating) random
boundaries studied in \cite{gomez2011wave}, and to electromagnetic
waves studied in
\cite{alonso2015electromagnetic,marcuse2013theory}. Waveguides with
slowly changing cross-section can be considered as well, especially if
the cross-section increases in the forward direction of wave
propagation.  Otherwise, the problem is complicated by the presence of
turning waves \cite{anyanwu1978asymptotic} analyzed in random
waveguides in \cite{borcea2017pulse,borcea2017transport}.

The goal of imaging is to locate a penetrable scatterer called ``the
target'', which lies inside the waveguide filled with a homogeneous
medium, with density $\rho_o$ and bulk modulus $K_o$. 
The pressure field $p(t,\bx)$ and the velocity field ${\itbf u}(t,\bx)$ 
satisfy the acoustic wave equations
\ba
\label{eq:wave1a}
\frac{1}{K(\bx)} \partial_t p(t,\bx)+\nabla\cdot {\itbf u}(t,\bx) =0 ,\\
\rho_o \partial_t {\itbf u}(t,\bx)+\nabla p(t,\bx)= {\itbf F}(t,\bx),
\label{eq:wave1b}
\ea
for time $t \in \RR$ and  location $\bx$ inside the waveguide.
For simplicity, we assume that the target has no contrast density but its
bulk modulus $K(\bx)$ is different from the background:
\ba
\frac{1}{K(\bx)} = \frac{1}{K_o}[ 1 + {R}(\bx) \big]. 
\ea
The function ${R}(\bx)$ models the reflectivity of the target.  We use
throughout the orthogonal system of coordinates $\bx = (x,z)$ shown in
Figure \ref{fig:setup}, with range $z \in \RR$ measured along the
direction of the axis of the waveguide, starting from the source, and
with cross-range $x$ in the interval $\big( \cD^-(z), \cD^+(z)\big)$,
called the cross-section of the waveguide. 
The source term is of the form
\ba
{\itbf F}(t,\bx) = {\itbf e}_z f(t,x) \delta(z) e^{-i \omega_o t}+c.c.,
\label{eq:force}
\ea
where $c.c.$ stands for complex conjugate and ${\itbf e}_z$ is the unit vector
pointing in the $z$-direction.
The target is assumed
infinitesimally thin, with reflectivity
\ba
{R}(\bx) = {r}(x) \delta(z-L),
\label{eq:Int2}
\ea
where ${r}(x)$ is compactly supported in $\big( \cD^-(L),
\cD^+(L)\big)$. This assumption is convenient for the analysis, but the
results extend to targets of finite range support.

Combining equations \eqref{eq:wave1a} and \eqref{eq:wave1b},  we obtain that the pressure wave field is the solution 
of the wave equation
\ba
\Delta p(t,\bx) - \frac{1}{c^2(\bx)} \partial_t^2 p(t,\bx) = e^{-i
  \om_o t}f(t,x) \delta'(z) +c.c.,
\label{eq:Int3}
\ea
for $t \in \RR$ and  $\bx \in \big( \cD^-(z), \cD^+(z)\big) \times
\RR$,  with Dirichlet boundary conditions
$$
p(t,(x,z))=0,\quad x\in \{ \cD^-(z), \cD^+(z) \}, \, z \in \RR.
$$
The   wave speed $c(\bx)$  is defined by 
\ba
\frac{1}{c^2(\bx)} = \frac{1}{c_o^2}[ 1 + {R}(\bx) \big],
\label{eq:Int1}
\ea
where $c_o =\sqrt{K_o/\rho_o}$  is the constant  wave speed
in the homogeneous medium that fills the waveguide. 

The cross-section of the waveguide is randomly fluctuating and has constant mean width $X$. We model the boundary by
\ba
\cD^+(z) = X\big[ 1 + \sigma^+\mu^+(z) \big], \quad \cD^-(z) = {X}
\sigma^- \mu^-(z),
\label{eq:Int4}
\ea
using the stationary random processes $\mu^\pm(z)$ with mean zero $
\EE[\mu^\pm(z)] = 0 $, {where $\EE[\cdot]$ denotes the expectation
with respect to the statistical distribution of the random boundary
fluctuations.} The processes $\mu^\pm(z)$ may be
independent or correlated and are mixing, with rapidly decaying mixing
rate as defined in \cite[Section 2]{papanicolaou1974asymptotic}. Their
autocorrelation
\ba
\cR^\pm(z) = \EE \big[\mu^\pm (z ) \mu^\pm(0)\big],
\label{eq:Int6}
\ea
is normalized to peak value $\cR^\pm(0) = 1$, and its integral 
\ba
\int_\RR {\rm d} z \, \cR^\pm(z) = O(\ell),
\label{eq:Int7}
\ea
defines the typical range scale of the fluctuations, the correlation
length $\ell$. We assume as in \cite{alonso2011wave} that $\mu^\pm(z)$
are twice continuously differentiable, with almost surely bounded
derivatives, and quantify the standard deviation of the fluctuations
of $\cD^\pm(z)$ using the small, positive and dimensionless parameters
$\sigma^\pm \ll 1$.

The wave is generated by the random source \eqref{eq:force}  localized at $z = 0$ and  oriented in the range
direction. The source signal is modulated at the carrier (central)
frequency $\om_o$ and has a slowly varying envelope $f(t,x)$. This is
a complex-valued, stationary in time Gaussian process, with zero
mean
%\footnote{We denote throughout by $\EE[\cdot]$ the expectation  with respect to the statistical distribution of the random boundary
 % fluctuations and by $\lb \cdot \rb$ the expectation with respect to   the distribution of the random source.}
\ba
\lb f(t,x) \rb = 0,
\label{eq:Int8}
\ea
zero relation function 
\ba
\lb f(t,x) f(t',x')\rb = 0,
\label{eq:Int9}
\ea
and covariance 
\ba
\lb f(t,x) \overline{f(t',x')}\rb = \sqrt{B}F\big[B(t-t')] \theta(x,x').
\label{eq:Int10}
\ea
The bar is used throughout the paper to denote complex conjugate
{and $\lb \cdot \rb$ denotes the expectation with respect to
  the distribution of the random source.}
The function
$F$ in the expression of the covariance \eqref{eq:Int10} is real valued, bounded, and integrable.
The frequency scale $B$
in its argument is called the bandwidth, because it determines the
support of the power spectral density function of $f(t,x)$, which is
proportional to the Fourier transform of $F$ \cite{grimmet}. The real
valued, integrable function $\theta(x,x')$ models the spatial
coherence of the source. It is supported in $\cAs \times \cAs$, where
the interval $\cAs \subseteq \big(\cD^-(0),\cD^+(0)\big)$ is called the
source aperture. If the source is incoherent, then
\ba
\theta(x,x') = \delta(x-x')1_{\cAs}(x), \quad ~~
1_{\cAs}(x) = \left\{ \begin{array}{ll}1, \quad &x \in \cAs, \\
0,  & x \notin \cAs.
\end{array} \right.
\label{eq:Int11}
\ea

The data for imaging the target are gathered at a detector located at
range $z = L + \cL$, with support (aperture) in the interval $\cAd
\subseteq \big(\cD^-(L+\cL), \cD^+(L+\cL)\big).$ In conventional
imaging, the detector would measure the wave field $p(t,x,L+\cL),$ for $t
\in \RR$ and $x \in \cAd$, and the image would be formed using matched
field \cite{baggeroer1993overview} or reverse time migration
processing \cite{borcea2013quantitative,borcea2015inverse}. Such
imaging is impeded by scattering at the random boundary of the
waveguide.  The long range cumulative effect of this scattering is
described mathematically by the randomization of the wave and the
energy exchange between its components, called the wave modes. This
energy exchange is quantified by transport equations derived in
\cite{alonso2011wave} for waveguides with random boundaries, and in
\cite[Chapter 20]{fouque07} for waveguides filled with random media.
Imaging based on these equations is studied in
\cite{borcea2010source,acosta2015source}.

None of the aforementioned imaging methods apply to the setting
considered here, where the detector-target range offset $\cL$ is so
large that cumulative scattering distributes the energy evenly among
the wave modes. It is not useful to measure $p(t,x,L+\cL)$ in this
equipartition regime, so we consider instead measurements of the 
net energy flux through the detector
\ba
I(t) =  \int_{\cAd}{\rm d} x \, p(t,x,L+\cL) {\itbf e}_z \cdot {\itbf u}(t,x,L+\cL), \quad t \in \RR.
\label{eq:Int12}
\ea
This flux, which is also called the sound power in the physical litterature \cite{landau},
 is associated with the usual energy density  \cite[Section 2.1.8]{fouque07}
\ban
e(t,x,z) = \frac{1}{2 K_o} p(t,x,z)^2 + \frac{\rho_o}{2} |{\itbf u}(t,x,z)|^2.
\ean
In spite of the limited data \eqref{eq:Int12} and the strong scattering, we show that
it is possible to image the target using a ghost-like imaging
modality.  Ghost imaging was introduced in the optics literature
\cite{pittman1995,cheng2009ghost,li2010ghost,shapiro2012physics} and was analyzed
in the context of wave propagation in random media in
\cite{garnier2016ghost}. The random source is realized in this context
by a laser beam passed through a rotating glass diffuser
\cite{shapiro2012physics}, followed by a beam splitter that divides
the beam in two parts: The first part illuminates the object of
interest, which is typically a mask, and is then captured by a single
pixel (bucket) detector. The second part does not interact with the
mask and its time and spatially resolved energy flux is measured at a
high resolution detector. The ghost image of the mask is then formed
by correlation of the two energy flux measurements.

Here we have a different setting, where the whole source beam
illuminates the target, which is not a mask, but a penetrable
scatterer. 
{The goal is to build an imaging modality for monitoring changes (i.e., emergence of targets) in the
waveguide at $z = L$. If more than one range is of interest, then an
image may be formed at each such range. } The image is formed using the correlation of the data
\eqref{eq:Int12} with
\vspace{0.05in}
\ba
\Ir(t,x) = 
 \pr(t,x,L) {\itbf e}_z \cdot {\itbf u}^{({\rm r})}(t,x,L), 
\label{eq:Int13}
\ea
the energy flux 
at range $L$, in a reference waveguide {(hence the superscript (r))}, for the same wave
source.   We consider two reference
waveguides:
\begin{enumerate}
\itemsep 0.05in
\item The unperturbed waveguide, with straight boundary at $x = 0$ and
  $x = {X}$, where we can calculate $\Ir(t,x)$ analytically,  {provided we know the wave generated by the source.
 This implies either controlling the acoustic source or measuring the wave near the source.} 
  \item The waveguide with the same random boundary \eqref{eq:Int4}, 
where $\Ir(t,x)$ must be measured before the presence of the target.
\end{enumerate}
\vspace{0.05in}
In either case, the imaging function at range $z = L$ is defined by
\ba
\cI_{\tau,T}(x) = \int_0^\tau {\rm d} s \, \cC_T(s,x), 
\label{eq:Int14}
\ea
where $\cC_T$ is the empirical energy flux correlation
\ba
\cC_T(s,x) &=& \frac{1}{T} \int_0^T {\rm d} t \,
\Ir(t,x)\big(I(t+s)-\Iin(t+s)\big) - \Big[\frac{1}{T} \int_0^T {\rm d} t
  \, \Ir(t,x)\Big] \nonumber \\ 
  &&\times \Big[\frac{1}{T} \int_0^T
  {\rm d} t \, \big(I(t)-\Iin(t)\big)\Big]
\label{eq:Int15}
\ea
and $\Iin(t)$ is the incident energy flux at the detector {(hence the superscript (i))}, in the
absence of the target. This can be measured, or the terms involving it
in \eqref{eq:Int15} can be calculated in terms of the second-order
statistics of the random processes $\mu^\pm(z)$. 

{The imaging function \eqref{eq:Int14} is designed to reflect the transverse profile of the 
target at $z=L$. It is not the classical ghost imaging function used in optics \cite{shapiro2012physics},
where the waves are monochromatic, there is no waveguide effect and one can simply evaluate 
 $\cC_T(s= 0,x) $. Waveguides are dispersive, meaning that different components (modes) of the wave field 
 propagate at different speed along the axis of the waveguide. 
Our analysis shows that due to mode dispersion,   the imaging function  should involve the integration of $\cC_T(s,x) $
over the time lag $s$,  in a time window of duration $\tau$ that is long enough to 
encompass the arrival of sufficiently many propagating modes.}

{More explicitly, we show that  the two user defined parameters $\tau$ and $T$ affect the statistical stability
and focusing of the imaging function \eqref{eq:Int14} at the target.  Statistical stability means that the
image is approximately the same for all the realizations of the random
source and boundary,
\ba
\cI_{\tau,T}(x) \approx \EE \big[ \lb \cI_{\tau,T}(x) \rb \big].
\label{eq:Int14_2}
\ea
We will see that the empirical correlation \eqref{eq:Int15} converges
as $T \to \infty$, in probability, to the statistical correlation with
respect to the distribution of the random source.  We will also see that for a large
enough bandwidth $B$, the imaging function is insensitive to the
realization of the fluctuations $\mu^\pm(z)$. 
Thus, to ensure a robust image, the wave field should not 
be monochromatic and the integration time $T$ should be large enough.
Once the bandwidth is taken into account, the mode dispersion in the 
waveguide  becomes important. The  time parameter $\tau$ controls how many 
modes are used in the image formation and consequently, it determines 
the resolution  of the image. At the very least, $\tau$ 
should satisfy the order relation
\ba
\frac{\cL}{c_o} < \tau \ll T,
\label{eq:Int14_1}
\ea
where the lower bound defines the travel time scale from the target to
the detector. 
If we want all the propagating modes to contribute, it should satisfy } 
\ba
\frac{\cL}{c_o} \ll \tau \ll T.
\label{eq:Int14_1b}
\ea

The goal of the paper is to analyze the imaging function
\eqref{eq:Int14}. The analysis is based on the theory of wave
propagation in random waveguides, developed in \cite{alonso2011wave}.
We quantify the resolution of the image in terms of the type of
reference waveguide, the target range $L$, the apertures $\cAs$ and
$\cAd$ of the source and detector, the partial coherence of the
source, and the time parameter $\tau$.

The paper is organized as follows: We begin in section \ref{sect:Form}
with the mathematical model of long range wave propagation in random
waveguides and the mathematical model of the measurements \eqref{eq:Int12} and the reference 
energy flux \eqref{eq:Int13}. The results of the analysis of the imaging function \eqref{eq:Int14} are in
section \ref{sect:Imag} and the details of the calculations are given in appendices.  
We end with a summary in section
\ref{sect:Sum}.
% ---------
\section{Long range wave propagation in the random waveguide}
\label{sect:Form}
In this section we summarize the wave propagation results obtained in
\cite{alonso2011wave}. As explained in
\cite{alonso2011wave,borcea2015inverse} these results are
qualitatively (but not quantitatively) the same as those in waveguides
filled with random media analyzed in
\cite{kohler77,dozier1978statistics} and \cite[Chapter 20]{fouque07}.
We begin with the scaling regime in section \ref{sect:scale}. The
mathematical model of the wave in the empty random waveguide, called
the incident wave, is given in section \ref{sect:incid}. The model of
the scattered wave is in section \ref{sect:scatt}.  We use it to derive the mathematical model 
of \eqref{eq:Int12}--\eqref{eq:Int13} in section \ref{sect:meas}.
\subsection{Long range scaling}
\label{sect:scale}
Let us introduce the small and positive dimensionless parameter $\ep
\ll 1$, which gives the order of magnitude of the standard deviation
of the random fluctuations of the boundary in \eqref{eq:Int4},
\ba
\sigma^\pm \sim \ep,
\label{eq:LR1}
\ea
where the symbol $\sim$ means that $\sigma^\pm/\ep$ is bounded above and
below by positive constants independent of $\ep$.

The effect of these fluctuations depends on the relation between the
correlation length $\ell$ and the carrier wavelength $\la_o = 2 \pi c_o/ \om_o$. We assume as in \cite{alonso2011wave} that
\ba
\ell \sim \la_o,
\label{eq:LR2}
\ea
so there is an 
efficient interaction of the wave with the random boundary.
Because the standard deviation \eqref{eq:LR1} is small,  this interaction 
 has a negligible effect at small range,
\ba
p(t,\bx) \approx p_o(t,\bx), \qquad z \sim \ell,
\label{eq:LR3}
\ea
where $p_o(t,\bx)$ is the solution of the wave equation in the
unperturbed waveguide. However, the scattering effect accumulates as
the wave propagates and becomes significant at $z \sim \ell/\ep^2$. The
wave $p(t,\bx)$ is randomized at such ranges and it is quite different
from $p_o(t,\bx)$, as shown in \cite{alonso2011wave}.

We denote by $z_\ep$ the range coordinate in this scaling, with the
target and detector located at $z_\ep = \Lep$ and $z_\ep = \Lep +
\cLep$, where
\ba
\frac{\ep^2 \Lep}{\ell} \sim 1, \quad \frac{\ep^2 \cLep}{\ell} \sim 1.
\label{eq:LR4}
\ea
In an abuse of notation, we let 
\ba
\label{eqLscaled}
\Lep = \frac{L}{\ep^2}, \quad \cLep = \frac{\cL}{\ep^2},
\ea
with $L$ and $\cL$ independent of $\ep$, satisfying 
\ba
\frac{\cL}{L} \gg 1,
\label{eq:LR5}
\ea
which means that the target-detector distance is much larger than the source-target distance.

The mean width of the cross-section  of the waveguide satisfies 
\ba
\frac{{X}}{\la_o} \gg 1,  \quad \mbox{independent of } \ep,
\label{eq:LR6}
\ea
so that the wave has many propagating components (modes). This is
needed to achieve a good resolution of the image, as long as the time
parameter $\tau$ in \eqref{eq:Int14} is chosen appropriately.  We
rename this parameter $\tau_\ep$, to emphasize its dependence on
$\ep$, and obtain from \eqref{eq:Int14_2} and
\eqref{eq:LR4}--\eqref{eq:LR5} that it should satisfy the scaling
relation
\ba
\tau_\ep = \frac{\cT}{\ep^2}, 
\label{eq:LR6p}
\ea
with $ \cT$ independent of $\ep$. The other time parameter
$T$, used in the calculation of the empirical correlation
\eqref{eq:Int15}, is much larger than $\tau_\ep$, so we can analyze
the imaging function by taking first the limit $T \to \infty$ and then
$\ep \to 0$.

To have the same number of propagating modes at all the frequencies in
the spectrum of the source, we take a small $\ep$ dependent bandwidth 
denoted  by $B_\ep$,  
\ba
\frac{B_\ep}{\om_o} \sim \ep^\alpha, \quad \alpha \in (1,2).
\label{eq:LR7}
\ea
The choice $\alpha < 2$ ensures that the covariance \eqref{eq:Int10}
is supported at time offsets $\sim 1/B_\ep$ that are much smaller than
the  travel time of order $\Lep/c_o$ from the source to the target. Therefore, we may
think of $F$ in \eqref{eq:Int10} as a pulse. The choice $\alpha < 2$ also gives
the statistical stability of the image \eqref{eq:Int14} with respect
to the realizations of $\mu^\pm(z_\ep)$, as explained in section
\ref{sect:statstab}.  The choice $\alpha > 1$ is convenient because we
can neglect some deterministic mode dispersion effects in the
waveguide.

\subsection{The incident wave}
\label{sect:incid}
As is usual in scattering theory, we call the solution
$\pin(t,x,z_\ep)$ of the wave equation in the empty waveguide ``the
incident wave'', {hence the index ${(i)}$}. We write it using the Green's function of the
Helmholtz equation
\ba
\big[\Delta + k^2(\om+\om_o)\big]\hat g\big(\om, (x,z_\ep),
(\xi,\zeta_\ep)\big) = \delta(x-\xi)\delta'(z_\ep-\zeta_\ep),  
\label{eq:IW0}
\ea
at frequency offset $\om$ from the carrier $\om_o$, where
$k(\om+\om_o) = (\om+\om_o)/c_o$ is the wavenumber.  This Green's
function is bounded and outgoing at $|z_\ep| \to \infty$ and defines
the incident wave as
\ba
\pin(t,x,z_\ep) &=& \frac{1}{2 \pi}\int_{-\infty}^\infty e^{-i
  (\om+\om_o)t} \, {\rm d}\hpin(\om,x,z_\ep) +c.c.,
\label{eq:IW9} \\
{\rm d}\hpin(\om,x,z_\ep) &=& \int_{\cAs} 
  \hat g\big(\om, (x,z_\ep), (\xi,0)\big) {\rm
  d} \hat f(\om,\xi) {\rm d} \xi   .
\label{eq:IW}
\ea
In this expression we used the spectral theorem for stationary
processes \cite[Section 9.4]{grimmet} to represent $f(t,x)$ by the
complex Gaussian measure ${\rm d} \hat f(\om,x)$ satisfying
\ba
\lb {\rm d} \hat f(\om,x) \rb &=& 0, \label{eq:IWf1}\\
 \lb {\rm d} \hat
f(\om,x) {\rm d} \hat f(\om',x') \rb &=& 0, \label{eq:IWf2}\\ 
\lb {\rm  d} \hat f(\om,x) {\rm d}\overline{ \hat f(\om',x')} \rb &=& \frac{2
  \pi}{\sqrt{B_\ep}} \hat F \Big(\frac{\om}{B_\ep}\Big)
\theta(x,x')\delta(\om-\om') {\rm d} \om{\rm d} \om'.
\label{eq:IWf3} 
\ea
We call ${\rm d} \hat f(\om,x)$ the Fourier transform of $f(t,x)$, in
an abuse of terminology.  Note that $\hat F$, the Fourier transform of
$F$, is the power spectral density of a stationary process, so it is
real valued and non-negative by Bochner's theorem \cite{grimmet}.

The Green's function is given by a superposition of time harmonic
waves (modes) which are either propagating or evanescent. The mode
decomposition is obtained by expanding $\hat g$ in the orthonormal
basis $\big(\phi_j(x)\big)_{j \ge 1}$ of the eigenfunctions
\ba
\phi_j(x) = \sqrt{\frac{2}{{X}}} \sin \Big(\frac{\pi j x}{{X}} \Big)
\label{eq:IW1}
\ea
of the operator $\partial_x^2$ with Dirichlet boundary conditions at
$x = 0$ and $x = X$. This expansion is justified because equation
\eqref{eq:IW0} can be rewritten in the unperturbed waveguide domain
$(0,{X})\times \RR$ with a change of variables that flattens the
boundary and maps the fluctuations $\mu^\pm(z_\ep)$ to coefficients in
the transformed Laplacian operator \cite{alonso2011wave}.  The
expression of the Green's function at $z_\ep > \zeta_\ep$ is
\ba
\hat g (\om,(x,z_\ep),(\xi,\zeta_\ep)) & \approx &\sum_{j=1}^{N(\om)}
\phi_j(x)\frac{a_{j}\big(\om,z_\ep;(\xi,
  \zeta_\ep)\big)}{\sqrt{\beta_j(\om + \om_o)}} e^{i
  \beta_j(\om+\om_o)(z_\ep-\zeta_\ep)} \nonumber \\
  &&+ \sum_{j =
  N(\om)+1}^\infty \phi_j(x) \hat g_{j}\big(\om,z_\ep; (\xi,
\zeta_\ep)\big),
\label{eq:IW2}
\ea
where the approximation is because we neglect the backward propagating
modes.

The first term in \eqref{eq:IW2} sums the $N(\om) $ forward
propagating modes. With our choice \eqref{eq:LR7} of the bandwidth we
have $\om \sim B_\ep \sim \ep^\alpha \om_o$ and therefore,
\ba
N(\om) = \left\lfloor \frac{k(\om + \om_o) {X}}{\pi}
\right\rfloor \approx N(\om_o) , \quad \mbox{for}~ \ep \ll 1,
\label{eq:IW3}
\ea
where $\lfloor \cdot \rfloor$ denotes the integer part. We omit
henceforth the argument $\om_o$ of $N$, to simplify notation.  The
$j$th forward propagating mode in \eqref{eq:IW2} is the combination
of two plane waves
\ba
 \phi_j(x)e^{i \beta_j(\om+\om_o)(z_\ep-\zeta_\ep)} = -
 \frac{i}{\sqrt{2{X}}} \left[ e^{ i \bka_j^+ \cdot (x,z_\ep-\zeta_\ep)} - e^{
     i \bka_j^- \cdot (x,z_\ep-\zeta_\ep)} \right],
\label{eq:IW4}
\ea
that cancel each other at the boundary points $x =0$ and $x = X$. Their wave vectors 
\ban
\bka_j^\pm = \big(\pm \frac{\pi j}{{X}}, \beta_j(\om+\om_o) \big)
\ean
have Euclidian norm $\|\bka_j^\pm\| = k(\om+\om_o)$ and a positive
component in the range direction, called the $j$th mode wavenumber,
\ba
\beta_j(\om+\om_o) = \sqrt{k^2(\om+\om_o) - \Big(\frac{\pi j}{{X}}\Big)^2 }.
\label{eq:IW6}
\ea
The backward propagating waves have a similar expression, except that
their wave vectors have range components of opposite sign. Because we
assume smooth random boundary fluctuations, we can use the forward
scattering approximation which neglects these backward going waves. A
detailed justification of this approximation is in
\cite{alonso2011wave} (see also \cite{kohler77} and \cite[Chapter
  20]{fouque07}).

The terms $\hat g_{j}$ of the series in \eqref{eq:IW2} are the
evanescent modes, which decay exponentially with the range offset
$z_\ep - \zeta_\ep$ and are negligible at the detector. However,
they interact with the propagating modes over the long range of
propagation, as described in \cite{alonso2011wave}. This interaction
is taken into account in the statistical description of the
propagating mode amplitudes $a_{j}$.  These are random fields with
starting values
\ba
a_{j}\big(\om,\zeta_\ep;(\xi,\zeta_\ep) \big) =
\frac{\sqrt{\beta_j(\om+\om_o)}}{2} \phi_j(\xi), \quad j = 1, \ldots, N, 
\label{eq:IW7}
\ea
and evolve at $z_\ep > \zeta_\ep$ as described by
\ba
\left( \begin{array}{c}
a_{1}\big(\om,z_\ep;(\xi,\zeta_\ep)\big) \\
\vdots \\
a_{N}\big(\om,z_\ep;(\xi,\zeta_\ep)\big)
\end{array}
\right)
= 
{\bf P}(\om,z_\ep;\zeta_\ep) \left( \begin{array}{c}
a_{1}\big(\om,\zeta_\ep;(\xi,\zeta_\ep)\big) \\ \vdots
\\
a_{N}\big(\om,\zeta_\ep;(\xi,\zeta_\ep)\big)
\end{array}
\right),
\label{eq:IW8}
\ea
using the complex, random $N \times N$ propagator matrix
${\bf P}(\om,z_\ep;\zeta_\ep)$.  This matrix equals the identity at
$z_\ep = \zeta_\ep$ and its statistics is described in the limit $\ep
\to 0$ in \cite{alonso2011wave}.

\subsection{The scattered wave}
\label{sect:scatt}
The solution 
of the wave equation \eqref{eq:Int3}, with wave speed defined in
\eqref{eq:Int1}--\eqref{eq:Int2}, has the form 
\ba
p(t,x,z_\ep) = \frac{1}{2 \pi}\int_{-\infty}^\infty e^{-i
  (\om+\om_o)t} \, {\rm d} \hat p(\om,x,z_\ep) +c.c.,
\label{eq:IW10}
\ea
where ${\rm d} \hat{p}$ satisfies the Lippmann-Schwinger equation
\ba
{\rm d} \hat p(\om,x,z_\ep) &=& {\rm d} \hpin(\om,x,z_\ep) -
k^2(\om+\om_o) \int_0^{X} {\rm d} y \, {r}(y) \hat
G\big(\om,(x,z_\ep),(y,L_\ep)\big) \nonumber \\
&& \times {\rm d}\hat
p(\om,y,L_\ep).
\label{eq:IW11}
\ea
Here we introduced another Green's function, satisfying the Helmholtz
equation
\ba
\big[\Delta + k^2(\om+\om_o)\big]\hat G\big(\om, (x,z_\ep),
(\xi,\zeta_\ep)\big) = \delta(x-\xi)\delta(z_\ep-\zeta_\ep).
\label{eq:IW12}
\ea
This equation is like \eqref{eq:IW0}, except that there is no range derivative
of the Dirac $\delta(z_\ep-\zeta_\ep)$. The expression of $\hat G$ is
similar to \eqref{eq:IW2},
\ba
\hat G\big(\om, (x,z_\ep), (\xi,\zeta_\ep)\big) \approx
\sum_{j=1}^{N} \phi_j(x)\frac{A_{j}\big(\om,z_\ep;(\xi,
  \zeta_\ep)\big)}{\sqrt{\beta_j(\om + \om_o)}} e^{i
  \beta_j(\om+\om_o)(z_\ep-\zeta_\ep)},
\label{eq:IW13}
\ea
where we neglected the evanescent modes. The random mode amplitudes
$A_j$ are defined as in equation \eqref{eq:IW8}, by the same
propagator ${\bf P}(\om,z_\ep;\zeta_\ep)$, but their starting values  are different,
\ba
A_j\big(\om,\zeta_\ep;(\xi, \zeta_\ep)\big) = \frac{\phi_j(\xi)}{2 i
  \sqrt{\beta_j(\om+\om_o)}}, \quad j = 1, \ldots, N.
\label{eq:IW14}
\ea

The scattered wave
\ba
\psc(t,x,z_\ep) = \frac{1}{2 \pi} \int_{-\infty}^\infty e^{-i
  (\om+\om_o)t} \, {\rm d}\hpsc (\om,x,z_\ep) +c.c., 
\label{eq:IW10p}
\ea
with Fourier transform 
\ba
{\rm d} \hpsc(\om,x,z_\ep) = {\rm d}\hat p(\om,x,z_\ep) - {\rm
  d}\hpin(\om,x,z_\ep),
\label{eq:IW10pp}
\ea
depends nonlinearly on the reflectivity ${r}(y)$ of the target,
assumed small. In imaging it is common to use the Born approximation
of this wave, obtained by replacing ${\rm d} \hat p(\om,y,L_\ep)$ with
${\rm d} \hpin(\om,y,L_\ep)$ in the integrand in \eqref{eq:IW11}. The
imaging function \eqref{eq:Int14} depends quadratically on ${r}(y)$,
so for consistency, we keep the second order terms in the Born series
expansion of the scattered wave
\ba
{\rm d} \hpsc(\om,x,z_\ep)  &\approx & - k^2(\om+\om_o) \int_0^{X}
 \hspace{-0.05in}  {r}(y) \hat G\big(\om
 ,(x,z_\ep),(y,L_\ep)\big){\rm d} \hpin(\om,y,L_\ep) {\rm d}y\nonumber \\ &&+
 k^4(\om+\om_o) \iint_0^{X}\hspace{-0.05in} 
 {r}(y) {r}(y')\hat G\big(\om,(x,z_\ep),(y,L_\ep)\big)  \nonumber \\ && \hspace{0.5in}
 \times \hat
 G_o\big(\om,(y,L_\ep),(y',L_\ep)\big) {\rm d}\hpin(\om,y',L_\ep){\rm d} y {\rm d} y' .
\label{eq:Born}
\ea
Here we replaced $\hat G\big(\om ,(y,L_\ep),(y',L_\ep)\big)$ by
the Green's function in the unperturbed waveguide
\ba
 \hat G_o\big(\om,(y,L_\ep),(y',L_\ep)\big) = \sum_{j=1}^\infty
 \frac{\phi_j(y)\phi_j(y')}{2 i \beta_j(\om + \om_o)},
\label{eq:Born1}
\ea
because the random boundary has a negligible effect on the wave
propagation between the nearby points $(y,L_\ep)$ and $(y',L_\ep)$ in the
support of the target. The
series in \eqref{eq:Born1} involves both the propagating and the
evanescent modes, with wavenumbers $\beta_j(\om+\om_o)$ defined in
\eqref{eq:IW6} for $j = 1, \ldots, N$, and with
\ba
\beta_j(\om+\om_o) = i \sqrt{\Big(\frac{\pi j}{{X}}\Big)^2 -
  k^2(\om+\om_o)}, \quad j > N.
\label{eq:IW6ev}
\ea

\subsection{The measurements and the reference energy flux}
\label{sect:meas}
The detector measures  the spatially integrated energy flux 
\ba
I(t) = \int_{\cAd} {\rm d}x \, p(t,x,L_\ep+\cL_\ep) {\itbf e}_z \cdot {\itbf u}(t,x,L_\ep+\cL_\ep) ,
\label{eq:Im1}
\ea
with $p(t,x,L_\ep+\cL_\ep)$ defined in equations \eqref{eq:IW10},
\eqref{eq:IW10p}--\eqref{eq:Born}. These measurements are correlated in 
\eqref{eq:Int14} with the energy flux 
\ba
\Ir(t,x) = \pr(t,x,L_\ep ) {\itbf e}_z \cdot {\itbf u}^{({\rm r})} (t,x,L_\ep )  ,
\label{eq:Im2}
\ea
of the wave, where
\ba
&& \hspace{-0.5in}\pr(t,x,L_\ep) = \frac{1}{2 \pi} \int_{-\infty}^\infty \int_{\cAs}  
e^{-i (\om + \om_o)t}   \hgr\big(\om, (x,L_\ep), (\xi,0)\big) {\rm d} \hat f(\om,\xi) {\rm d} \xi ,
\ea 
and 
\ba
 && \hspace{-0.5in} {\itbf e}_z \cdot {\itbf u}^{({\rm r})} (t,x,L_\ep)=  \frac{1}{2 \pi} \int_{-\infty}^\infty \int_{\cAs}  
\frac{e^{-i (\om + \om_o)t}}{i (\omega_o+\omega)} \partial_z \hgr\big(\om, (x,L_\ep), (\xi,0)\big) {\rm d} \hat f(\om,\xi) {\rm d} \xi   .
\label{eq:Im3}
\ea
The expression of $\pr(t,x,L_\ep)$  is similar to \eqref{eq:IW}, except for the Green's
function $\hgr$ which models the wave propagation in the reference
waveguide.
The expression of $ {\itbf e}_z \cdot {\itbf u}^{({\rm r})} (t,x,L_\ep)$ follows from (\ref{eq:wave1b}) 
which reads in the Fourier domain as
\ban
-i \rho_o (\omega+\omega_o) {\itbf e}_z \cdot {\rm d} \hat{\itbf u}^{({\rm r})}(\om,x,z_\ep)  + \partial_z {\rm d} \hat{p}^{({\rm r})}(\om,x,z_\ep) = 0 .
\ean

We consider two reference waveguides: The  unperturbed
waveguide, where
\ba 
\hgr\big(\om , (x,L_\ep), (\xi,0)\big)  \approx \frac{1}{2} \sum_{j =
  1}^N \phi_j(x) \phi_j(\xi) e^{i \beta_j(\om+\om_o)L_\ep},
\label{eq:Im4}
\ea
and the  random, empty waveguide, where 
\ba 
&& \hspace{-0.5in} \hgr\big(\om , (x,L_\ep), (\xi,0)\big) = 
\hat g\big(\om , (x,L_\ep), (\xi,0)\big) \nonumber \\
&&  \hspace{0.93in}\approx \frac{1}{2} \sum_{j,q =
  1}^N \sqrt{\frac{\beta_q(\om_o)}{\beta_j(\om_o)}}\phi_j(x)
\phi_q(\xi) {P}_{jq}(\om ,L_\ep;0)e^{i \beta_j(\om+\om_o)L_\ep}.
\label{eq:Im5}
\ea
In both cases we neglect the evanescent modes, which have a very small
contribution at the large range $L_\ep$.  The random mode amplitudes
are written in \eqref{eq:Im5}  in terms of the entries ${P}_{jq}$ of
the propagator, as in equations \eqref{eq:IW7}--\eqref{eq:IW8}.
Because  $\om/\om_o \sim \ep^\alpha$, with $\alpha \in
(1,2)$, we approximate the mode wavenumbers in the amplitudes by their
value at the carrier frequency. The phase in
\eqref{eq:Im4}--\eqref{eq:Im5} is
\ba
\beta_j(\om+\om_o) L_\ep = \beta_j(\om_o) L_\ep + \om \beta'_j(\om_o)
L_\ep + O\big(\ep^{2(\alpha - 1)}\big),
\label{eq:Im6}
\ea
where 
\ba
\frac{1}{\beta_j'(\om_o)} = \Big(\frac{{\rm d} \beta_j(\om+\om_o)}{{\rm d}\omega}\mid_{\omega=0}\Big)^{-1}= \frac{c_o \beta_j(\om_o)}{k_o},
\label{eq:Im7}
\ea 
is the $j$th mode speed, satisfying 
\ba
0<\frac{1}{\beta_N'(\om_o)} < \frac{1}{\beta_N'(\om_o)} < \ldots <
\frac{1}{\beta_1'(\om_o)} < c_o.
\label{eq:Im7p}
\ea 
We use throughout the notation $k_o = k(\om_o) = \om_o/c_o$. 

\section{Analysis of the ghost imaging function}
\label{sect:Imag}
The imaging function \eqref{eq:Int14} is defined by the time integral
of the empirical correlation of the reference energy flux $\Ir(t,x)$ and
the difference $I(t)- \Iin(t)$ of the net intensities at the
detector. The subtraction of the net energy flux $\Iin(t)$ of the
incident wave is so that the imaging function vanishes in the absence
of the target. This energy flux could be measured in the empty
waveguide or, if this is not feasible, its contribution to the
imaging function can be estimated in terms of the second order
statistics of the random boundary fluctuations $\mu^\pm(z_\ep)$, as
explained in section \ref{sect:Incid}.

The empirical energy flux correlation \eqref{eq:Int15} converges in the
limit $T \to \infty$, in probability, to its statistical expectation
with respect to the distribution of the random source
\cite[Proposition 2.3]{garnier2016passive}. Assuming a large 
integration time $T$, we obtain the approximation
\ba
&& \hspace*{-1.in}
\cC_T(s,x) \approx \cC(s,x) = \lb
\Ir(t,x)\big[I(t+s)-\Iin(t+s)\big]\rb - \lb \Ir(t,x)\rb \lb
I(t)-\Iin(t)\rb.
\label{eq:Im8}
\ea
This expression is independent of $t$, because of the time
stationarity of the source. We calculate it in 
\ref{app:A}.

The imaging function is given by 
\ba
\cI_{\tau_\ep,T}(x) \approx \cI_{\tau_\ep}(x) = \int_0^{\tau_\ep} {\rm d} s
\, \cC(s,x),
\label{eq:Im10p}
\ea
with $\tau_\ep = \cT/\ep^2$ scaled as in \eqref{eq:LR6p}. Its expression
\ba
&&\hspace*{-1.in}
\cI_{\frac{\cT}{\ep^2}}(x)  
\approx  \frac{2 \cT}{(2 \pi)^2 \rho_o^2 c_o^2
  B_\ep}
  {\rm Re} \int_{-\infty}^\infty \hspace{-0.05in} {\rm d}\om \, \hat
F\Big(\frac{\om}{B_\ep}\Big)^2  \int_{-\infty}^\infty {\rm d}h \,
\mbox{sinc} \left(\frac{h \cT}{2}\right)e^{-i h \frac{\cT}{2}} \nonumber \\ 
&&\hspace*{-0.5in}\times 
     \iint_0^{X} \hspace{-0.05in}
     {\rm d} y {\rm d} y' \Big[ \cG_{12}(\om-\ep^2 h,y',x)
     \cG_{11}(\om,y,x) + \overline{\cG_{12}(\om,y,x)}
     \overline{\cG_{11}(\om-\ep^2 h,y',x)} \Big]
     \nonumber\\
     &&\hspace*{-0.5in}\times 
     \Big[ 
     -2  r(y) \cG_{31}( \omega,\omega-\ep^2 h,y,y')  
     -2 r(y') \cG_{32}( \omega,\omega-\ep^2 h,y,y')  
     \nonumber \\
     &&\hspace*{-0.5in}
     \quad +2 k_o^2 r(y) \int_0^{X} {\rm d}y''  r(y'') \cG_{31}( \omega,\omega-\ep^2 h,y'',y')  \hat{G}_o\big(\omega,(y'',L_\ep),(y,L_\ep)\big) 
          \nonumber\\
     && \hspace*{-0.5in}  
      \quad     +2 k_o^2 r(y') \int_0^{X}{\rm d}y''  r(y'') \cG_{32}( \omega,\omega-\ep^2 h,y,y'') \overline{ \hat{G}_o\big(\omega-\ep^2 h,(y'',L_\ep),(y',L_\ep)\big)}
     \nonumber\\
     &&\hspace*{-0.5in} \quad +   k_o^2 r(y) r(y') \cG_{2}( \omega,\omega-\ep^2 h,y,y') 
     \Big]
\label{eq:Im10}
\ea
is derived in  \ref{app:A}, where 
we recall that $\hat F$ and
${r}$ are real valued. Here we  introduced the
notation
\ba
&&\hspace*{-1.in}
\cG_{11}(\om,y,x) = \iint_{\cAs}\hspace{-0.05in} {\rm d} \xi {\rm d}
\xi' \, \theta(\xi,\xi') \hat g \big(\om,(y,L_\ep),(\xi,0)\big)
\overline{\hat{g}^{({\rm r})} \big(\om,(x,L_\ep),(\xi',0)\big)},
\label{eq:defM11}
\\
&&\hspace*{-1.in}
\cG_{12}(\om,y',x) = \iint_{\cAs}\hspace{-0.05in} {\rm d} \xi {\rm d}
\xi' \, \theta(\xi,\xi')\overline{ \hat g \big(\om,(y',L_\ep),(\xi,0)\big)}
(-i)\partial_z {\hat{g}^{({\rm r})} \big(\om,(x,L_\ep),(\xi',0)\big)},
\label{eq:defM12}
\\
&&\hspace*{-1.in}
\cG_2(\om,\om',y,y') = \int_{\cAd} \hspace{-0.06in}{\rm d} x' \, \hat
G\big(\om,(x',\Lep+\cLep),(y,\Lep)\big)  i \partial_z \overline{\hat
  G\big(\om',(x',\Lep+\cLep),(y',\Lep)\big)},
\label{eq:defM2}
\\
&&\hspace*{-1.in}
\cG_{31}(\om,\om',y,y') = \hspace{-0.03in}\int_{\cAd}\hspace{-0.06in}
       {\rm d} x' \, \hat G\big(\om,(x',\Lep+\cLep),(y,\Lep)\big)
      i \partial_z \overline{\hat g\big(\om',(x',\Lep+\cLep),(y',\Lep)\big)},
\label{eq:defM31}
\ea
and 
\ba
&&\hspace*{-1.in}
\cG_{32}(\om,\om',y,y') = \hspace{-0.03in}\int_{\cAd}\hspace{-0.06in}
       {\rm d} x' \, \hat g\big(\om,(x',\Lep+\cLep),(y,\Lep)\big)
      i \partial_z \overline{\hat G\big(\om',(x',\Lep+\cLep),(y',\Lep)\big)}.
\label{eq:defM32}
\ea

We calculate next the imaging function \eqref{eq:Im10} in both the
unperturbed and random waveguide, in order to quantify its focusing in
terms of the target range $\Lep$, the apertures $\cAs$ and $\cAd$, the
source coherence function $\theta$ and the time parameter $\cT$.

\subsection{Imaging in the unperturbed waveguide}
\label{sect:ImagH}
The expression of the imaging function is derived in  \ref{app:B}. It is 
given by the formula \eqref{eq:ImHom} 
for an arbitrary integration time $\cT$.  This $\cT$ does not play an important role 
in the support of the image, although it affects its magnitude at the target, but it must be larger than $\beta'_1(\om_o) \cL$,
the scaled travel time of the fastest 
propagating mode from the target to the receiver. Otherwise, $ \cI_{\frac{\cT}{\ep^2}}(x)\equiv 0$. 
In section \ref{sect:imHom}  we give the expression of the imaging function  for $
\cT > \beta'_{N}(\om_o) \cL,
$
which is simpler and independent of $\cT$, and  analyze in section \ref{sect:resHom} its focusing in terms of the 
source coherence and the source and detector apertures. The general function \eqref{eq:ImHom} is 
displayed in section \ref{sect:displayHom} for the case of a point target.
\subsubsection{The imaging function for a long integration time}
\label{sect:imHom}
\begin{figure}[t]
\begin{center}
\includegraphics[width=2.4in]{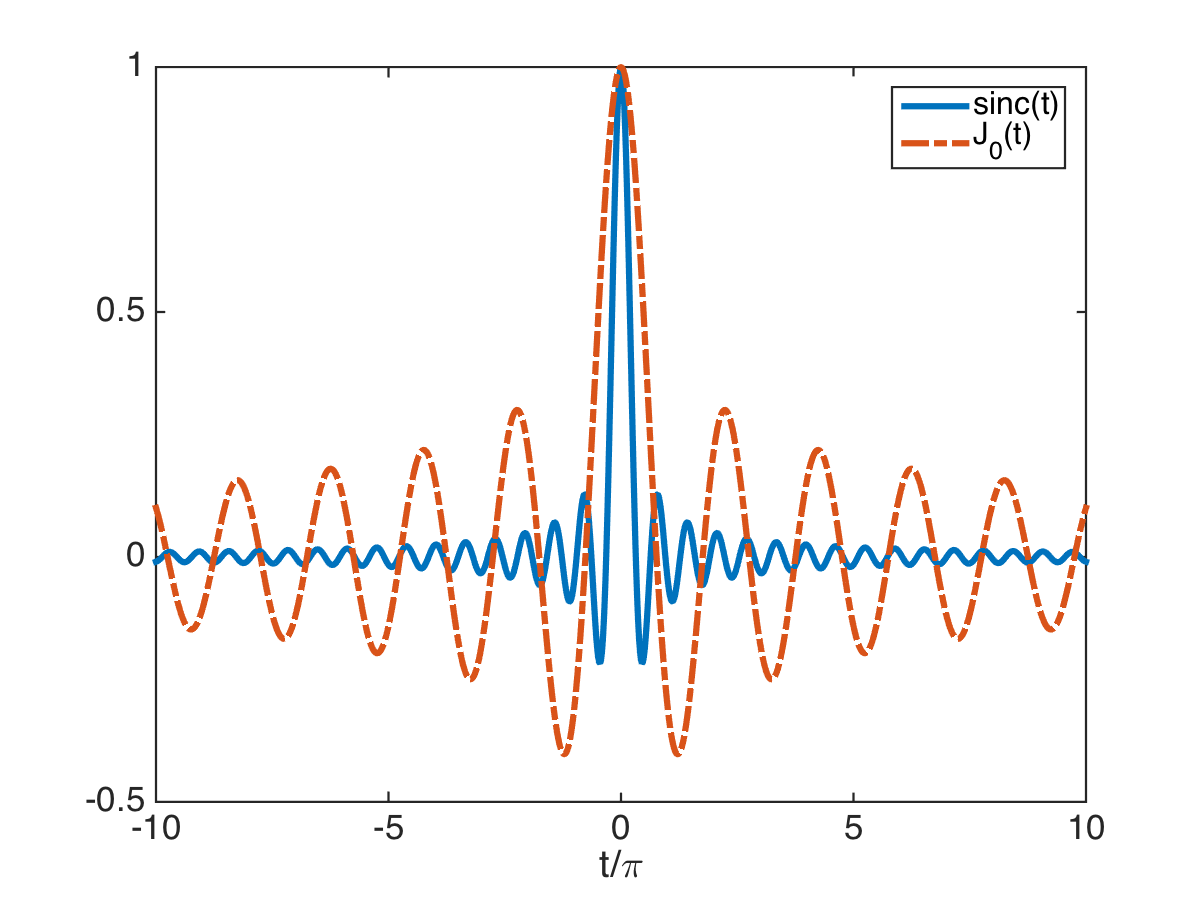}
\includegraphics[width=2.4in]{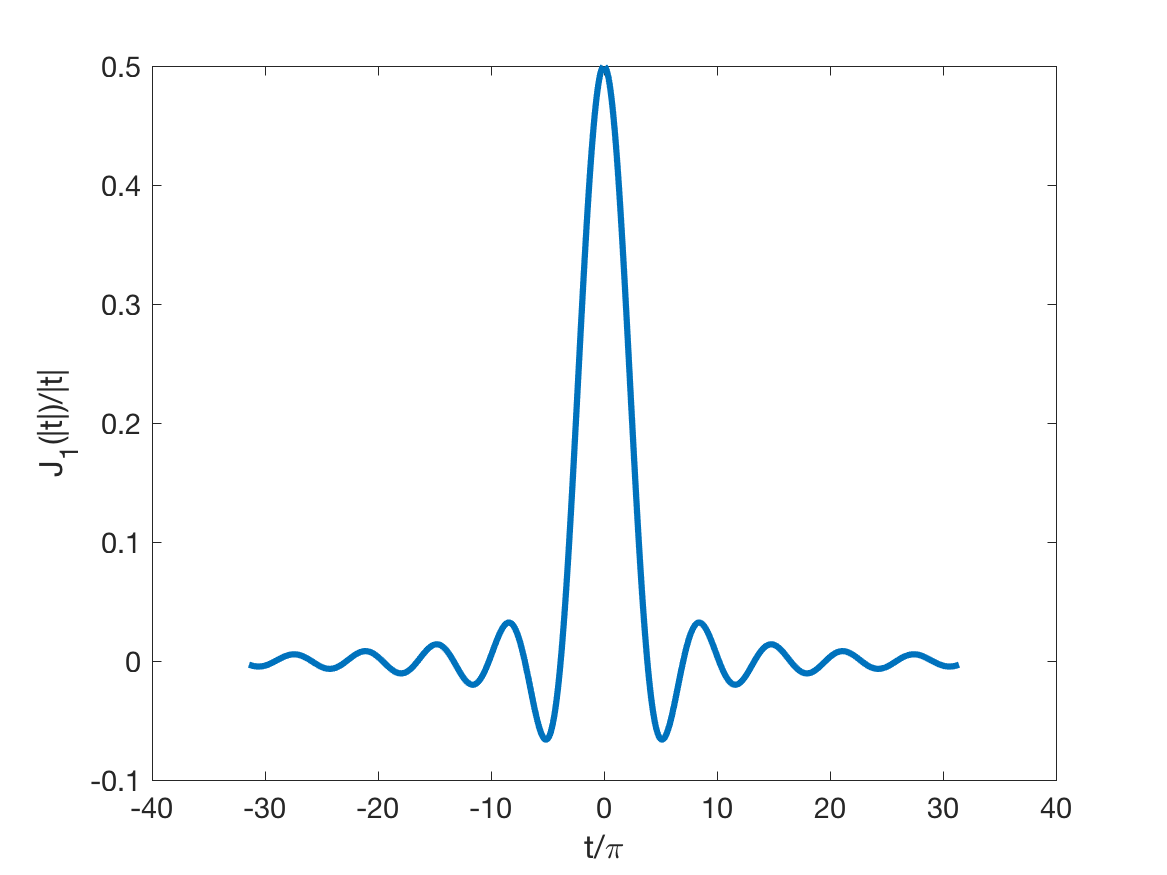}
\end{center} 
\caption{{The Bessel functions which approximate $\Phi_j$, $j=-1,0,1$, as in (\ref{eq:phi-1}-\ref{eq:phi0})}.
Left plot: the sinc function $\mbox{sinc}(t)$ (solid line) and the Bessel
  function $J_0(t)$ (dashed line) for $t \in [-10 \pi, 10 \pi]$. Right plot: the function $\frac{J_1(|t|)}{|t|}$ 
  for $t \in [-10 \pi, 10 \pi]$
  The argument $t$ is scaled by
  $\pi$ in the abscissa.}
\label{fig:sinc}
\end{figure}

The expression of the imaging function at $\cT > \beta'_{N}(\om_o) \cL$ is 
\ba
&&\hspace*{-1.in}
\cI (x) \approx  \frac{k^2_o \|F\|^2}{32 \rho_o^2c_o^2}
\iint_0^{X} {\rm d} y {\rm d}y'  \Big[
  \Phi_{{\bf S}0}(x,y) \Phi_{{\bf S}1}(x,y')+  \Phi_{{\bf S}1}(x,y) \Phi_{{\bf S}0}(x,y') 
\nonumber \\ 
\nonumber
&& \hspace*{1.1in} +2 \sum_{\stackrel{j,j'=1}{j\neq j'}}^N 
  \beta_j(\om_o) {S}_{jj'}^2 \phi_j(x)^2 \phi_{j'}(y) \phi_{j'}(y')  \Big]\\
 &&\hspace*{-0.4in}
 \times
\Big[ -2 \int_0^{X} {\rm d} y'' \, \Phi_{{\bf D}0}(y',y'') \Phi_{-1}(y,y'') {r}(y) r(y'')
+\Phi_{{\bf D}-1}(y,y') {r}(y){r}(y')\Big],
  \label{eq:Im11}
\ea
where $\|F\|$  is the $L^2$ norm
of $F$ and  ${\bf S}$ and ${\bf D}$  are the 
$N \times N$ matrices  with
entries
\ba
&&\hspace*{-1.in}
{S}_{jj'} = \iint_{\cAs} {\rm d} \xi {\rm d}\xi' \, \theta(\xi,\xi')
\phi_j(\xi) \phi_{j'}(\xi'), \qquad  {D}_{jj'} = \int_{\cAd}
    {\rm d} x \, \phi_j(x) \phi_{j'}(x), \label{eq:H3}
\ea
that depend on the source and receiver apertures and the coherence
function $\theta$ of the source.  Here we introduced the sums $\Phi_n$, 
$\Phi_{{\bf S} n}$, $\Phi_{{\bf D} n}$ defined  for $n \geq -1$ by
\ba
 \Phi_{ n}(x,y) &=&\sum_{j=1}^N  \beta_j(\om_o)^n
 \phi_j(x) \phi_j(y) ,
 \label{eq:Im13}\\
 \label{eq:PhiS}
\Phi_{{\bf S} n}(x,y) &=& \sum_{j=1}^N \beta_j(\om_o)^n {S}_{jj}
  \phi_j(x) \phi_j(y), \\
\Phi_{{\bf D} n}(x,y) &=&\sum_{j=1}^N \beta_j(\om_o)^n {D}_{jj}
  \phi_j(x) \phi_j(y).
\label{eq:PhiD}
\ea
Note that when the source is incoherent (delta-correlated) with $\theta$ defined in \eqref{eq:Int11}, and has 
full aperture $\cAs = (0,X)$, then $\bS$ is the identity matrix. Similarly, $\bD$ equals the identity when 
$\cAd = (0,X)$. In these ideal cases the sums \eqref{eq:PhiS}--\eqref{eq:PhiD} equal \eqref{eq:Im13},
and can be approximated using that $N \gg 1$ by the scaling assumption
\eqref{eq:LR6}. Except very close to the boundaries $x,y \in \{0,X\}$ we have 
\ba
\Phi_{-1}(x,y) &
  \approx& \frac{2}{\pi N}
\sum_{l=1}^N \Big[1-\Big(\frac{l}{N}\Big)^2\Big]^{-\frac{1}{2}} \sin \Big(k_o x\frac{l}{N} \Big)\sin
  \Big(k_o y\frac{l}{N} \Big)
\nonumber \\ 
    &\approx & \frac{1}{2} 
  J_0\big(k_o(x-y)\big)  ,
  \label{eq:phi-1}
\ea
and similarly,
\ba
\Phi_0(x,y) \approx \frac{k_o}{\pi}   {\rm sinc}  ( k_o|x-y|)   , \qquad 
\Phi_1(x,y) \approx \frac{k_o^2}{4} \frac{J_1( k_o|x-y|)}{k_o|x-y|} , 
\label{eq:phi0}
%,\\  \Phi_2(x,y) =& \frac{2k_o^3}{\pi} \Big[ \frac{\sin s}{s^3} - \frac{\cos(s)}{s^2}\Big]_{s= k_o|x-y|}.
\ea
where $J_0$ and $J_1$ are the Bessel function of the first kind and of orders $0$ and $1$. All these functions are peaked 
at $x = y$, as illustrated in Fig. \ref{fig:sinc}.

\subsubsection{Quantification of resolution}
\label{sect:resHom}
The expression \eqref{eq:Im11} shows that the imaging function is  independent
of the target range $\Lep = L/\ep^2$.

The aperture $\cAd$ of the detector, which defines the matrix ${\bf D}$,  
has a marginal role. For example, in the case $\cAd = (0,a_d)$, 
with $a_d \le {X}$,
\ban
{D}_{jj} = \int_0^{a_d} {\rm d} x' \, \phi_j^2(x') = \frac{a_d}{{X}} \Big[1 - \mbox{sinc} \Big(\frac{2 j \pi a_d}{{X}}\Big)\Big] 
\approx \frac{a_d}{{X}}, 
\ean
where the approximation is for $j > {X}/a_d$.
Thus, the functions $\Phi_{{\bf D} n}$ in \eqref{eq:Im11} are approximately proportional to $\Phi_{n}$, for $n = -1, 0$, 
and they give a similar contribution to the focus of the image for a full detector aperture or 
a smaller one.

The focusing of \eqref{eq:Im11}  at points $x$ in the support of the reflectivity ${r}$ 
is primarily dictated by the coupling matrix ${\bf S}$.
The best focusing is for ${S}_{jj'} = \delta_{jj'}$, corresponding to
an incoherent source with $\theta$ defined in \eqref{eq:Int11} and a
full aperture $\cAs = (0,{X})$. The terms in  $\sum_{j\neq j'}$  in the square brackets in 
\eqref{eq:Im11} vanish in this case and, if in addition the detector has full aperture, the imaging function becomes 
\ba
\cI (x) 
\approx 
 -\frac{k^2_o  \|F\|^2}{16  \rho_o^2c_o^2}
\iint_0^{X} {\rm d} y {\rm d}y'  \, r(y) r(y') 
  \Phi_{0}(x,y) \Phi_{1}(x,y')\Phi_{-1}(y,y') . \label{eq:ImHomFULL}
\ea
The functions  $\Phi_n$, for $n = -1, 0,1$, are peaked at $x = y$ and are large for
\ban
|x-y| \le \frac{\pi}{k_o} = \frac{\la_o}{2},
\ean
as seen from \eqref{eq:phi-1}--\eqref{eq:phi0} and Fig. \ref{fig:sinc}.
Therefore,  the
image \eqref{eq:Im11} focuses at $x$ in the support of the
reflectivity ${r}$, with resolution $\la_o/2$.
%More explicitly, 
%if the reflector is point-like, i.e. ${r}(x)=\delta_{x_o}(x)$, then
%\ba
%\cI (x)  
%\approx& - \frac{k^5_o \|F\|^2}{256 \pi \rho_o^2c_o^2} \Psi( k_o |x-x_o|)  ,
%\ea
%where 
%\ba
%\Psi(s)=& \frac{2}{\pi} \int_{-\infty}^\infty {\rm sinc}(s') \frac{J_1(|s-s'|)}{|s-s'|} ds' {\rm sinc}(s)
%+\frac{2}{\pi}\int_{-\infty}^\infty {\rm sinc}(s') {\rm sinc}(|s-s'|) ds' \frac{J_1(s)}{s}
%\\
%\nonumber
%& - 2 {\rm sinc}(s)  \frac{J_1(s)}{s} \\
%=&2 {\rm sinc}(s)  \frac{J_1(s)}{s} ,
%\label{def:Psi}
%\ea
%where $\Psi(s)$ is an even function that peaks at $0$, with $\Psi(0)=1$.
Note that $\cI (x) $ has the form of a negative peak,
which means that the target appears as a shadow.

If the source does not have full aperture or if it is partially coherent,
then the imaging function remains focused at points in the support of
${r}$, as long as the matrix ${\bf S}$ is diagonally dominant.  
Otherwise, the terms in $\sum_{j\neq j'}$  in the first square brackets in
\eqref{eq:Im11}, which are not focused at points in the support of
${r}$, become significant.

\subsubsection{{Illustration of the point spread function}}
 \label{sect:displayHom}
 {Here we illustrate the resolution analysis given above by displaying in Figure \ref{fig:Res1} the 
imaging function $\cI(x)$ for a point target at cross-range $0.39X$, in a waveguide that supports $40$ propagating modes. 
We consider an incoherent source modeled as in \eqref{eq:Int11}, with $\cAs = (0,a_s)$, and a detector with aperture $\cAd = (0,a_d)$, for $0 < a_s, a_d  \le X$. }

 {To illustrate the effect of the duration $\cT$ of the integration window, we display in Figure \ref{fig:Res1} the image 
 $\cI(x)$ calculated using the general formula \eqref{eq:ImHom}, for three values of $\cT$: The first satisfies 
 $\cT > \beta_N' \cL$, as assumed in the previous section, and the result is shown with the solid blue line. 
 The other two satisfy $\beta_{n}' \cL < \cT <  \beta_{n+1}' \cL$ for $n = 10$ and $n = 1$, and the results are shown with the dotted red line and the solid black line. We also illustrate the effect of the aperture of the source and detector by displaying the images for 
$a_d = X$ and $a_s = X, 0.5X, 0.1X$ in the top plots and for $a_d = 0.5X, 0.1X$ and $a_s = X$ in the bottom plots.}

{\begin{figure}[t]
\begin{center}
\raisebox{0.05in}{\includegraphics[width=2.17in]{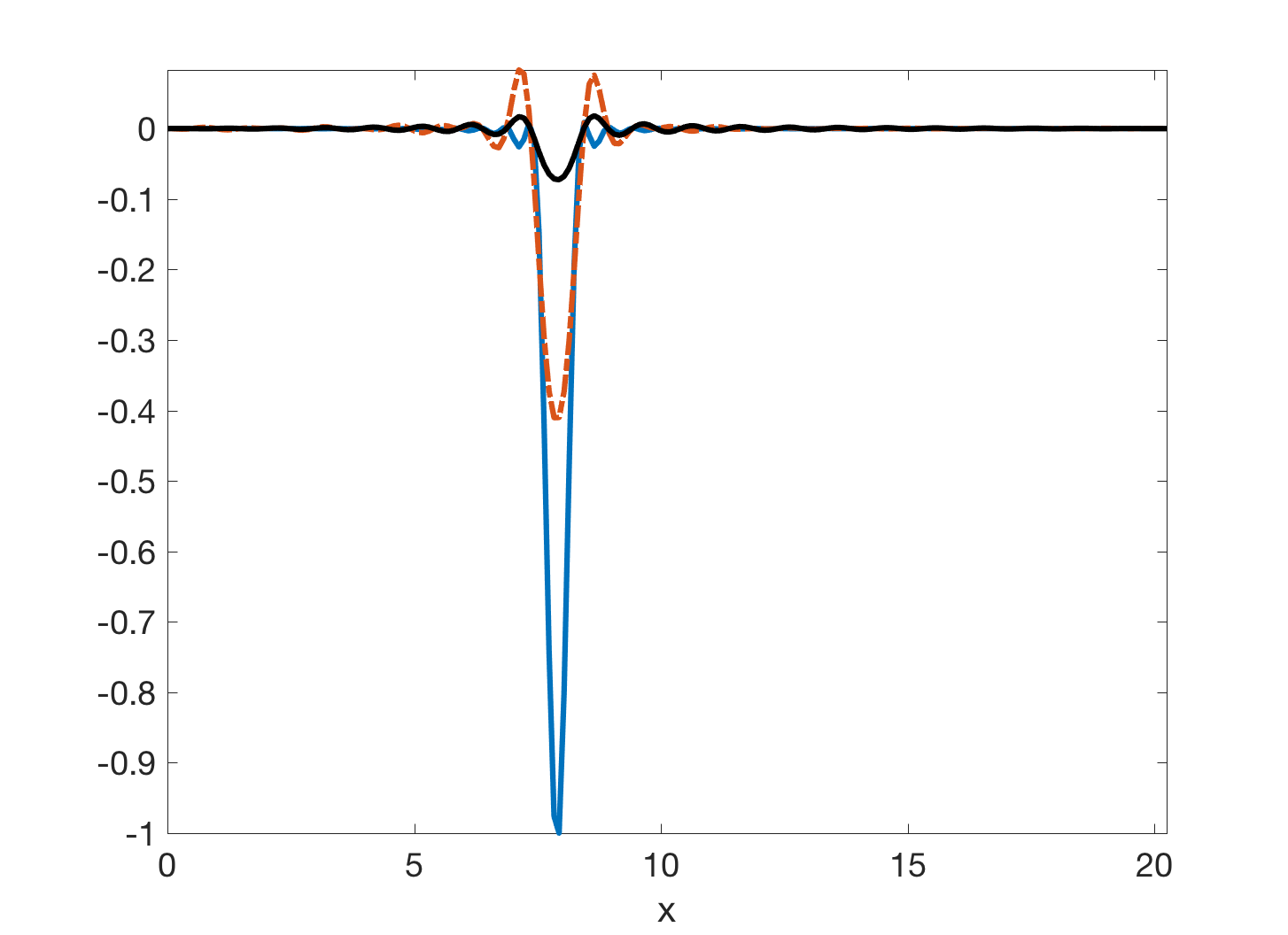}}
\hspace{-0.3in}
\includegraphics[width=2.17in]{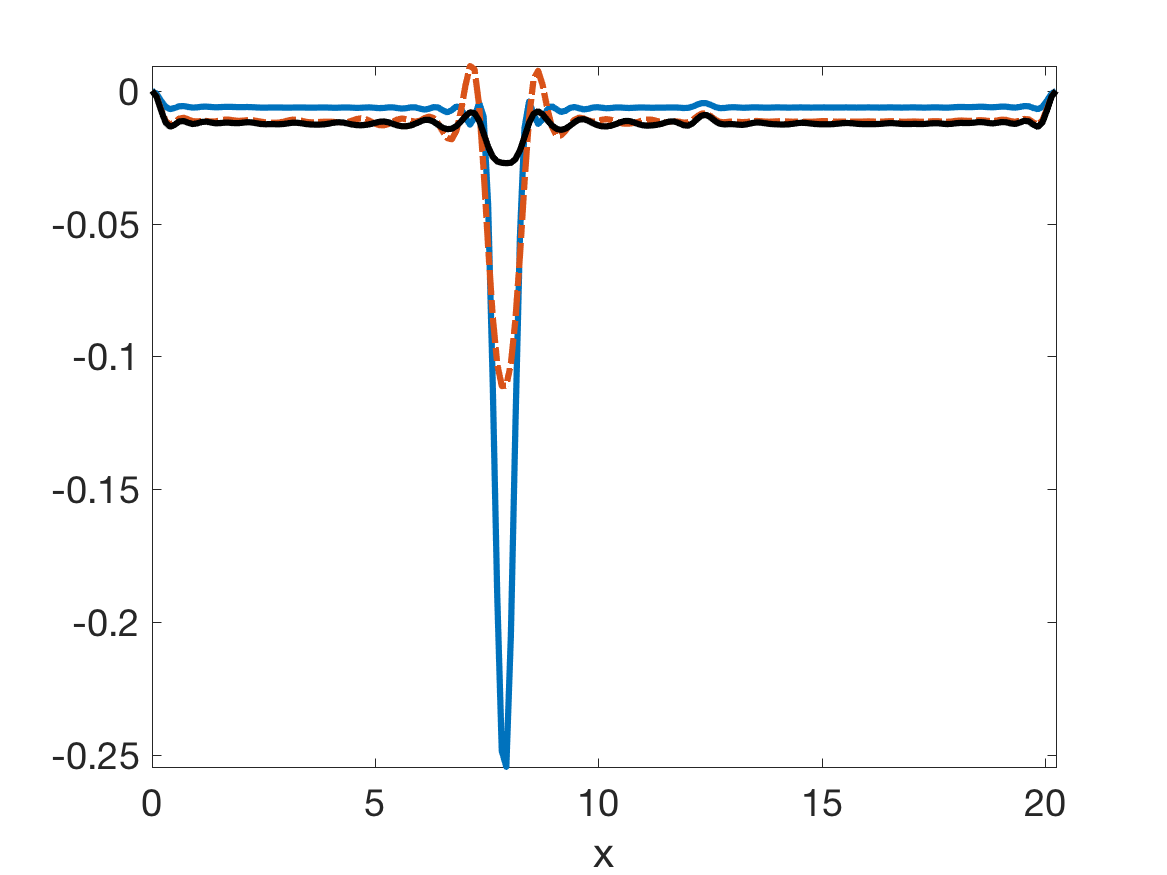}
\hspace{-0.3in}
\includegraphics[width=2.17in]{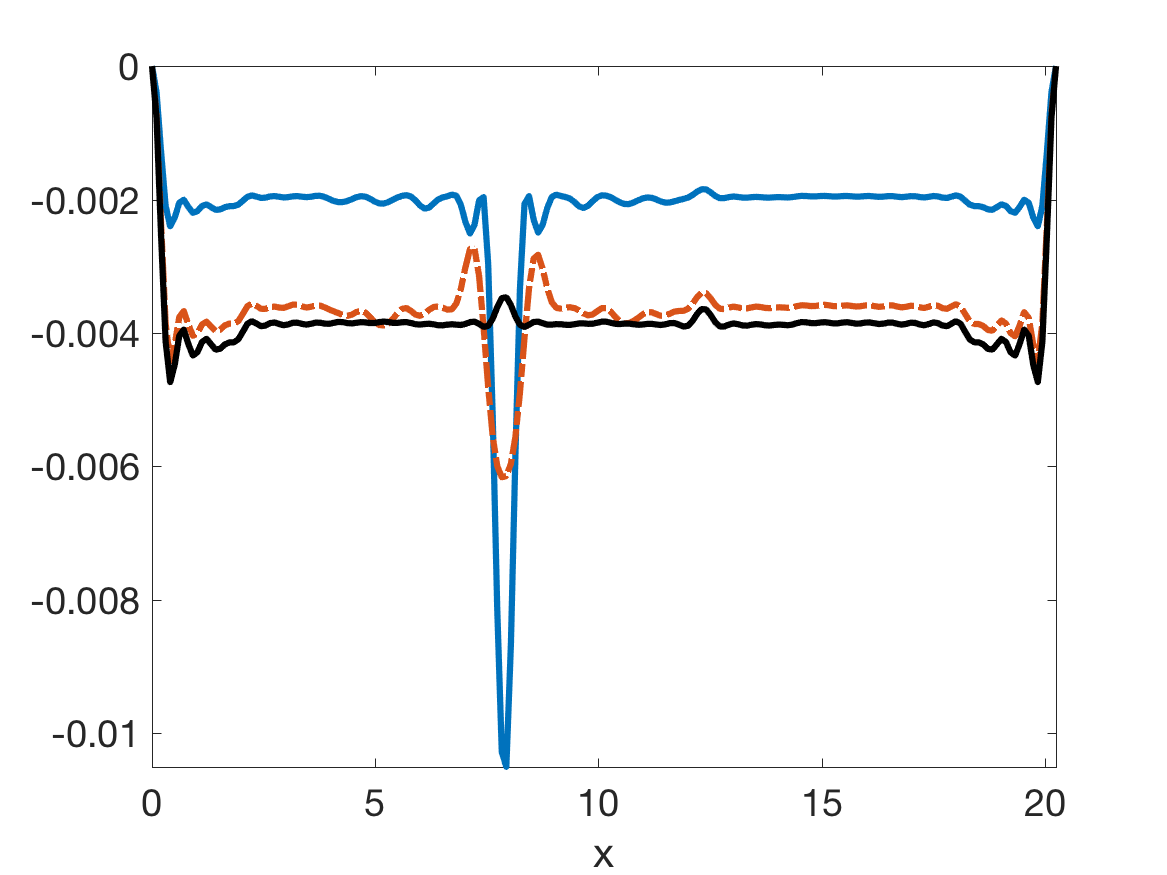} \\
\hspace{1.95in}\includegraphics[width=2.17in]{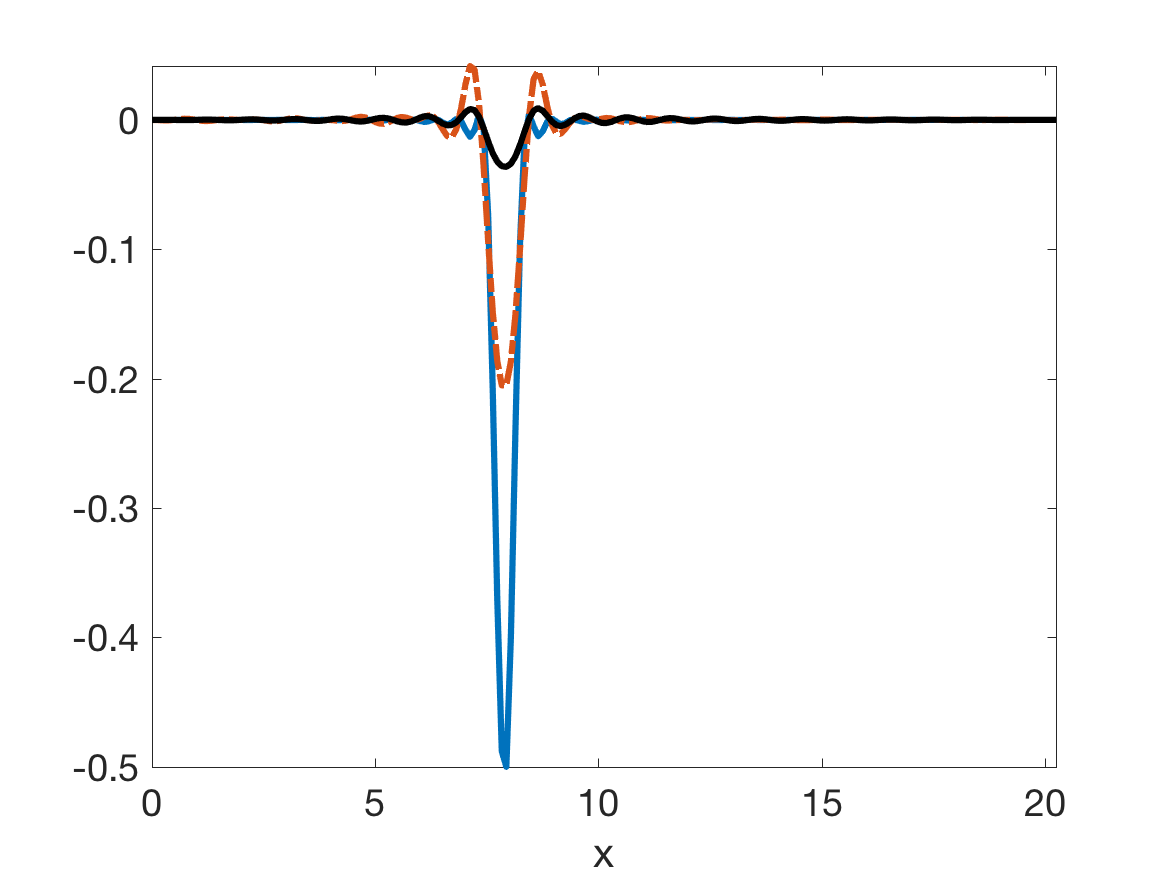}
\hspace{-0.22in}\includegraphics[width=2.17in]{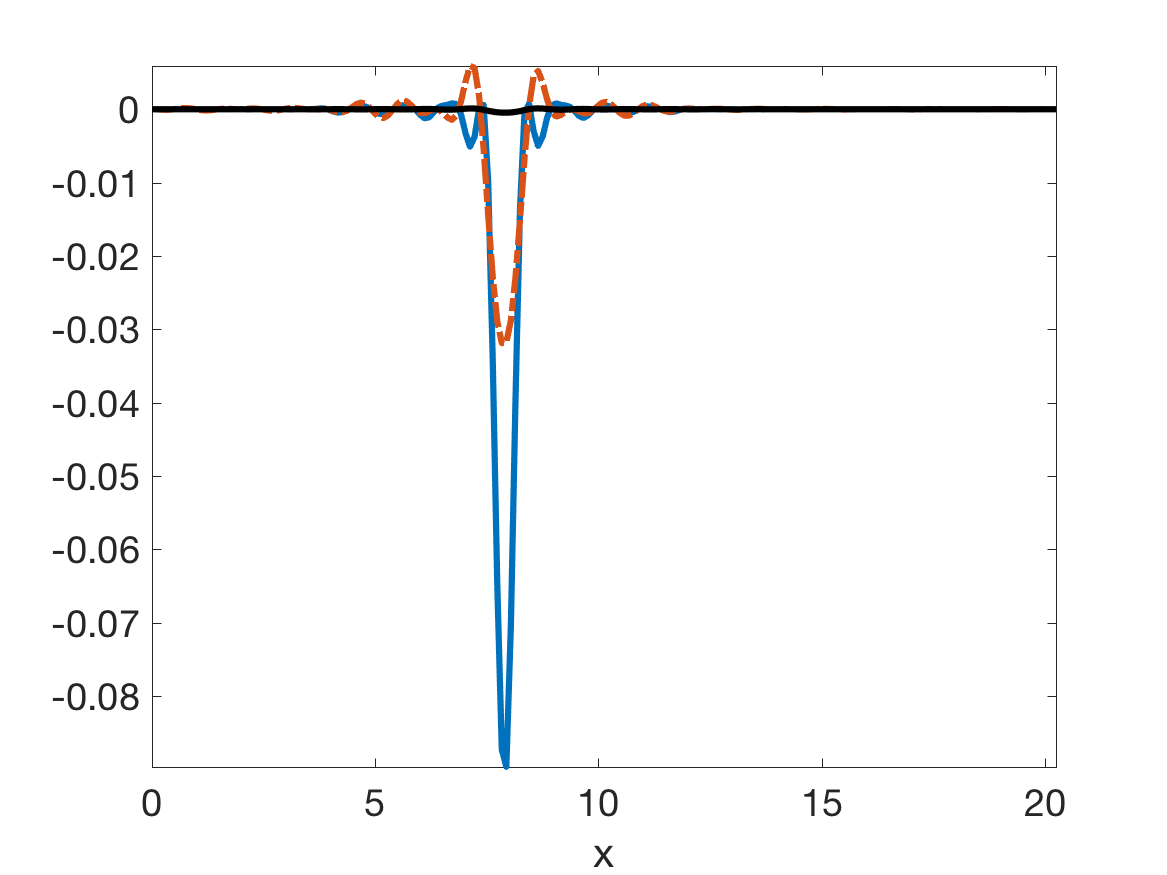}
\vspace{-0.1in}\caption{{
Imaging function \eqref{eq:ImHom}  calculated for $\cT > \beta'_N \cL$  \eqref{eq:Im11} (blue solid line) and 
$\beta_{n}' \cL < \cT <  \beta_{n+1}' \cL$ for $n = 10$ (red dotted line) and $n = 1$ (black solid line).
The top plots are  for $a_d = X$ and, from left to right, $a_s = X, 0.5X, 0.1X$.
The bottom plots are for $a_s = X$ and, from left to right, $a_d = 0.5X, 0.1X$.
We normalize by the absolute value at the target location 
of the image calculated with  $a_d = a_s = X$ and $\cT > \beta'_N \cL$. The abscissa 
is cross-range in units of the carrier wavelength.
}
}
\end{center}
\label{fig:Res1}
\end{figure}
}

{As explained in the previous section,  the image has a negative peak i.e., the target appears as a shadow 
at the true cross-range location $0.39X$. The value of $\cT$ has little effect on the support of this peak, but it 
affects its magnitude. The best results are for the large $\cT$, where the peak  is  more prominent and the side lobes are 
smaller. }

{The marginal effect of the detector aperture $a_d$ discussed in the previous section is seen in Figure \ref{fig:Res1}
to consist of a multiplicative factor that makes the peak less prominent for the smaller $a_d$.  However, the source aperture $a_s$ 
has a stronger effect on the image, which is no longer  supported in the vicinity of the target location. At $50\%$ aperture (top middle plot)
the value of the imaging function away from the target is small, but at $10\%$ aperture this value increases when compared to the peak and the 
target becomes less visible, especially for the smaller $\cT$.}

\subsection{Imaging in the random waveguide}
\label{sect:ImagMR}
We present here the analysis of the imaging function \eqref{eq:Im10}, in the
waveguide with random boundary, using the unperturbed reference
waveguide.  To calculate
$\cI_{\frac{\cT}{\ep^2}}(x)$, we recall from
\cite{alonso2011wave,borcea2015inverse} (see also \cite[Chapter
  20]{fouque07}) the relevant facts about the statistics of the wave
propagator ${\bf P}(\om,z_\ep;\zeta_\ep)$ in the limit $\ep \to 0$:
\vspace{0.05in}
\begin{enumerate}
\itemsep 0.05in
\item The propagators ${\bf P}(\om,z_\ep;\zeta_\ep)$ and
  ${\bf P}(\om,z'_\ep;\zeta'_\ep)$ for any two non-overlapping range
  intervals $(z_\ep,\zeta_\ep)$ and $(z'_\ep,\zeta_\ep')$ are
  uncorrelated.
\item The propagators ${\bf P}(\om,z_\ep;\zeta_\ep)$ and
  ${\bf P}(\om',z_\ep;\zeta_\ep)$ for any two frequencies $\om$ and $\om'$
  satisfying $|\om-\om'| \gg \ep^2 \om_o$ are uncorrelated.
  Moreover,  
  \ba
  {\bf P}(\om,z_\ep;\zeta_\ep)\approx {\bf P}(\om',z_\ep;\zeta_\ep), \quad  \forall \om, \om' \mbox{ satisfying } |\om-\om'| \ll \ep^2 \om_o.
  \ea
\item The propagator satisfies the second moment formula
\ba
\EE \Big[ {P}_{jq}(\om,z_\ep;\zeta_\ep) \overline{{P}_{j'q'}(\om-\ep^2
    h,z_\ep;\zeta_\ep)}\Big] \to \delta_{jq} \delta_{j'q'}
e^{\kappa_{jj'}(z-\zeta)} \nonumber \\
\hspace{0.5in}+ \delta_{jj'}\delta_{qq'}
\int_{-\infty}^\infty {\rm d} t \, W_j^{(q)}(\om_o,t,z-\zeta) e^{i h
  \big[t - \beta_j'(\om_o)(z-\zeta)\big]},
\label{eq:MomProp}
\ea
where we let $z_\ep = z/\ep^2$, $\zeta_\ep = \zeta/\ep^2$ and used the
symbol $\to$ to denote convergence in the limit $\ep \to 0$. 
\end{enumerate}

\vspace{0.05in}
\noindent The first term in the right-hand side of \eqref{eq:MomProp} models the
energy of the coherent part of the wave. Its evolution in range is described by
the complex matrix $(\kappa_{jj'})$ given explicitly in
\cite{alonso2011wave} and \cite[Section 6.2]{borcea2015inverse}, in
terms of the covariance of the random processes
$\mu^\pm$. This matrix satisfies
\ba
 \mbox{Re}[\kappa_{jj'}] < 0, \quad \kappa_{jj'} =
 \overline{\kappa_{j'j}}, \quad \forall j, j' = 1, \ldots, N,
\label{eq:MomProp1}
\ea
so the coherent term in \eqref{eq:MomProp} decays exponentially in
range.  The length scales $1/|\kappa_{jj}|$ of decay are called the scattering mean 
free paths of the modes. These scales are analyzed in \cite{alonso2011wave} and decrease monotonically with the 
mode index $j$.
The imaginary part of $(\kappa_{jj'})$ vanishes on the
diagonal, but it is non-zero otherwise. It accounts for dispersive
effects due to scattering at the random boundary.  The second term in
\eqref{eq:MomProp} models the incoherent part of the energy, defined
by the real valued, continuous density $W_j^{(q)}$ of the Wigner
transform described in \cite[Section 6.2]{borcea2015inverse} and
\cite[Proposition 20.7]{fouque07}.

{We consider a target at very long range offset from the detector, satisfying 
\ba
\label{eq:scaleLeq}
\cL \gg L_{\rm eq} \ge \frac{1}{|\kappa_{11}|},
\ea
where $L_{\rm eq}$ is the equipartition distance defined in
\cite[Section 20.3.3]{fouque07}. This is interesting because the wave reaching the detector is incoherent,  
meaning explicitly that 
\begin{equation}
\hspace{-1in}\EE\big[ \hat G(\om,(x,L_\ep + \cL_\ep),(y,L_\ep) \big] \approx 0, ~~
\EE\big[ \hat g(\om,(x,L_\ep + \cL_\ep),(y,L_\ep) \big] \approx 0, ~~ \forall x , y \in (0,X).
\label{eq:ADD1}
\end{equation}
Moreover,  the strong scattering in the range interval $z \in (L_\ep, L_\ep + \cL_\ep)$ 
distributes the energy evenly between the modes, independent of their amplitudes at the target.
This makes conventional imaging impossible, but as we  explain in section \ref{sect:randRes}, 
the target can still be located with the ghost imaging modality.}

\subsubsection{Statistical stability} \label{sect:statstab} 
In \eqref{eq:Im10} we integrate over the
frequency $\om$, which is of the order of the bandwidth $B_\ep$,
scaled as in \eqref{eq:LR7}. The propagator decorrelates over
 frequency offsets of order $ \ep^2 \om_o$ (item 2. above), so we integrate over
many frequency decorrelation intervals and obtain by the law of large
numbers the statistical stability result
\ba
\cI_{\frac{\cT}{\ep^2}}(x) \approx \EE \Big[ \cI_{\frac{\cT}{\ep^2}}(x) \Big].
\label{eq:ImR11}
\ea

{Note  that the integration over $s \in (0,\tau_\ep)$ in \eqref{eq:Im10}, with $\tau_\ep = \cT/\ep^2$, ensures that 
the imaging function is given by the superposition of products of two Green's functions (one with a complex conjugate) at frequencies 
$\om$ and $\om - \ep^2 h$, with $h \sim \om_o$.  The  result \eqref{eq:ImR11} holds because the expectation in the right hand side 
is not small.  If we did not integrate over $s$ in \eqref{eq:Im10}, or we had a much smaller $\tau_\ep$, then the 
integrand would consist of products of uncorrelated Green's functions, at frequencies offset by more than order $\ep^2 \om_o$. The expectation of such products is given by the 
product of the expectations of the Green's functions, which are negligible as in \eqref{eq:ADD1}.  The expectation of the terms in the integrand in \eqref{eq:Im10} involve the second moments of the wave,  which are not small,  
and this is why we obtain the result \eqref{eq:ImR11}.}

When calculating the expectation of \eqref{eq:Im10} we note that
$(\cG_{1n})_{n=1,2}$ depend on the range section $(0,\Lep)$ of the waveguide,
whereas $\cG_2$ and $(\cG_{3n})_{n=1,2}$ depend on the range section
$(\Lep,\Lep+\cL_\ep)$. By item 1. above, $\cG_{11}$ and $\cG_{12}$ are uncorrelated from
$\cG_2$,  $\cG_{31}$, and $\cG_{32}$, so to approximate \eqref{eq:Im10} we need the
expectations of \eqref{eq:defM2}--\eqref{eq:defM32} and
$\EE\Big[\cG_{11}(\om,y,x) {\cG_{12}(\om-\ep^2 h,y',x)}\Big]$.

We  consider large integration times $\cT \gg \cL/c_o$ so that, 
 as in the homogeneous case, the integral in $h$ in \eqref{eq:Im10}
becomes concentrated at $h=0$ and the imaging function becomes independent of $\cT$. As a consequence, we get from (\ref{eq:Im10}):
\ba
&&\hspace*{-1.in}
\cI (x)   \approx   \frac{1 }{ \pi \rho_o^2 c_o^2
  B_\ep}
  {\rm Re} \int_{-\infty}^\infty \hspace{-0.05in} {\rm d}\om \, \hat
F\Big(\frac{\om}{B_\ep}\Big)^2
     \iint_0^{X} \hspace{-0.05in}
     {\rm d} y {\rm d} y'  \nonumber\\
     &&\times \EE \Big[ \cG_{12}(\omega,y',x)
     \cG_{11}(\om,y,x) + \overline{\cG_{12}(\omega,y,x)}
     \overline{\cG_{11}(\omega,y',x)} \Big]
     \nonumber\\
     &&\times 
     \Big[ 
     -2  r(y) \EE[ \cG_{31}( \omega ,\omega,y,y')]  
     -2 r(y') \EE[\cG_{32}( \omega,\omega,y,y')  ]
     \nonumber \\
     &&
     \quad +2 k_o^2 r(y) \int_0^{X}{\rm d} y'' r(y'') \EE[ \cG_{31}( \omega,\omega,y'',y')]  \hat{G}_o\big(\omega,(y'',L_\ep),(y,L_\ep)\big)
          \nonumber\\
     &&   
      \quad     +2 k_o^2 r(y') \int_0^{X} {\rm d} y'' r(y'') \EE[ \cG_{32}( \omega,\omega,y,y'')]
       \overline{ \hat{G}_o\big(\omega,(y'',L_\ep),(y',L_\ep)\big)}
     \nonumber\\
     && \quad +   k_o^2 r(y) r(y') \EE[ \cG_{2}( \omega,\omega,y,y') ]
     \Big]   .
     \label{eq:meancI1}
\ea

\subsubsection{The imaging function} 
\label{sect:randRes}

The expression of the imaging function is given in \eqref{eq:ImagRand} for an 
arbitrary scaled range offset $L$ between the source and target. Here we write it 
in the two extreme cases, where the formulas are simpler:
\\
1) For weak scattering i.e., $L$ {smaller} than the scattering mean free path 
of all the propagating modes, the imaging function has the  form
\ba
&&\hspace*{-1.in}
\cI (x) 
\approx  \frac{k^2_o \|F\|^2 C_{{\bf D}} }{32 \rho_o^2c_o^2}
\iint_0^{X} {\rm d} y {\rm d}y'  \Big[
  \Phi_{{\bf S}0}(x,y) \Phi_{{\bf S}1}(x,y')+  \Phi_{{\bf S}1}(x,y) \Phi_{{\bf S}0}(x,y') 
\nonumber \\ 
\nonumber
&& \hspace*{0.4in}  +2 \sum_{j\neq j'} 
  \beta_j(\om_o) {S}_{jj'}^2 \phi_j(x)^2 \phi_{j'}(y) \phi_{j'}(y')  \Big]\\
&&\hspace{-0.35in}\times
\Big[ -2 r(y) \hspace{-0.05in}\int_0^{X} {\rm d} y'' r(y'') \Phi_{0}(y',y'') \Phi_{-1}(y,y'') 
+\Phi_{-1}(y,y') {r}(y){r}(y')\Big]  ,
\label{eq:expressimrand1}
\ea
where 
\ba
C_{{\bf D}} = \frac{1}{N} \sum_{m=1}^N {D}_{mm} = \frac{1}{N} \int_{\cAd}\Phi_0(x',x'){\rm d} x' .
\label{eq:Im13p}
\ea
\\
2) When scattering is strong between the source and the target i.e., $L$ is {larger} than the scattering mean free 
path of all the propagating modes,  the imaging function has the form
\ba
&&\hspace*{-1.in} \cI (x) 
 \approx    \frac{k^2_o \|F\|^2 C_{{\bf D}} }{16 \rho_o^2c_o^2}
\iint_0^{X} {\rm d} y {\rm d}y'  \Big[
 \sum_{j, j',l=1}^N 
  {S}_{j'l}^2  \phi_l(x)^2 \phi_j(y) \phi_{j}(y') \frac{\beta_{j'}(\om_o) \beta_{l}(\om_o)}{\beta_j(\om_o)}
\nonumber \\ 
&&\hspace{1.35in}\times  \int_{-\infty}^\infty {\rm d} t \, W_j^{(j')}(\om_o,t,L) \Big] \nonumber \\
&& \hspace{-0.4in} \times
\Big[ -2 r(y) \hspace{-0.05in}\int_0^{X} {\rm d} y'' r(y'')  \Phi_{0}(y',y'') \Phi_{-1}(y,y'') 
+\Phi_{-1}(y,y') {r}(y){r}(y')\Big] .
\label{eq:expressimrand2}
\ea
Moreover, when $L$ exceeds the equipartition distance, the Wigner transform becomes independent of the indexes 
$j,j'=1,\ldots, N$ and satisfies 
\ban
\int_{-\infty}^\infty {\rm d} t \, W_j^{(j')}(\om_o,t,L) \approx  \frac{1}{N}, \quad L \gg L_{\rm eq}.
\ean

As we explain below, the source aperture and coherence determine the focusing of the image 
in the weak scattering regime. In the strong scattering regime, the image is not focused 
at the target, even in the ideal case of an incoherent source with full  aperture. Intermediate 
scattering regimes are an interpolation of these two extreme cases and the image 
focuses at the target when $L$ is smaller than the scattering mean free paths of sufficiently many 
modes.

%
%\ba
%\nonumber
%\cI (x) 
%\approx  & \frac{k^2_o \|F\|^2 C_{{\bf D}} }{16 \rho_o^2c_o^2 N}
%\iint_0^{X} {\rm d} y {\rm d}y'  \Big[
% \sum_{j, j',l=1}^N 
%  {S}_{j'l}^2  \phi_l(x)^2 \phi_j(y) \phi_{j}(y') \frac{\beta_{j'}(\om_o) \beta_{l}(\om_o)}{\beta_j(\om_o)}
%   \Big]\\
%  &\times
%\Big[ -2 \int_0^{X} \Phi_{0}(y',y'') \Phi_{-1}(y,y'') {r}(y''){\rm d} y'' {r}(y)
%+\Phi_{-1}(y,y') {r}(y){r}(y')\Big] .
%\ea

\subsubsection{Quantification of resolution}
Note that the size of the detector appears only in the  constant amplitude factor $C_{{\bf D}}$, so the focusing of the 
image is entirely dependent on the source and the target range $L$.

In the weak scattering regime {(between the source and the target)}, 
the best result is obtained for an incoherent (delta-correlated) source  that spans the entire cross section
of the waveguide.  The matrix $\bS$ equals the identity in this case and   \eqref{eq:expressimrand1} becomes 
\ba
\cI(x) \approx - \frac{k_o^2 \|F\|^2 C_{\bD}}{16 \rho_o^2 c_o^2} \iint_0^X {\rm d}y {\rm d} y' r(y) r(y') \Phi_0(x,y) \Phi_1(x,y') \Phi_{-1}(y,y').
\label{eq:expressimrand3}
\ea
This is proportional to the imaging function \eqref{eq:ImHomFULL} in the homogeneous waveguide, at full detector aperture. 
This shows that the strong wave scattering between the target and the detector is beneficial
{and removes the marginal effect of a limited 
detector aperture $\cAd$ on the focusing of the image that we had observed in a homogeneous waveguide in sections 
\ref{sect:resHom}--\ref{sect:displayHom}.  The image has a negative peak and  has a resolution of order $\la_o/2$, as illustrated in the top left plot of Figure \ref{fig:Res1}.}

%The form of the source is important, as we need both delta-correlated sources that span the cross section and weak scattering
%between $0$ and $L$ to get an image of the target.
%Indeed, if the source covariance function is of the form (\ref{eq:Int11}), then
%the imaging function has the form
%\ba
%\nonumber
%\cI (x) 
%\approx & \frac{k^2_o \|F\|^2 C_{{\bf D}} }{32 \rho_o^2c_o^2}
%\iint_0^{X} {\rm d} y {\rm d}y'  \Big[
%  \Phi_{0}(x,y) \Phi_{1}(x,y')+  \Phi_{1}(x,y) \Phi_{0}(x,y') 
%   \Big]\\
%  &\times
%\Big[ -2 \int_0^{X} \Phi_{0}(y',y'') \Phi_{-1}(y,y'') {r}(y''){\rm d} y'' {r}(y)
%+\Phi_{-1}(y,y') {r}(y){r}(y')\Big] ,
%\ea
%when scattering is weak between the source and the target.
%If  the target is point-like, ${r}(x)=\delta_{x_o}(x)$, then one finds
%\ba
%\cI (x) 
%\approx - \frac{k^5_o \|F\|^2 C_{{\bf D}}}{256 \pi \rho_o^2c_o^2} \Psi( k_o |x-x_o|)   ,
%\ea
%where $\Psi$ is given by (\ref{def:Psi}). We then get a nice image of the shadow of the target.

In the strong scattering regime {between the source and the target},  the image does not focus, no matter what the source aperture and coherence are.
More precisely, in the ideal case of a delta-correlated source with full aperture, and for $L > L_{\rm eq}$, 
the expression \eqref{eq:expressimrand2} becomes 
\ba
\cI(x) \approx - \frac{k_o^2 \|F\|^2 C_{\bD}}{16 N\rho_o^2 c_o^2} \Phi_2(x,x) \iint_0^X {\rm d}y {\rm d} y' r(y) r(y')  \Phi_{-1}(y,y')^2.
\label{eq:expressimrand4}
\ea
This is approximately constant for $x$ away from the boundary points $0$ and $X$, {as  illustrated in the top left plot in Figure \ref{fig:Res2}.}

\subsection{Contribution of the incident energy flux to the image}
\label{sect:Incid}
In case that $\Iin(t)$ cannot be measured, we show here that 
it is possible to estimate its contribution to the imaging function
\ba
\cIin_{\frac{\cT}{\ep^2}}(x) &= \int_0^{\cT/\ep^2} {\rm d} s \Big[\lb 
I^{({\rm r})}(t,x) \Iin(t+s)\rb - \lb 
I^{({\rm r})}(t,x)\rb \lb \Iin(t+s)\rb\Big].
\label{eq:CI1}
\ea
We calculate this expression in  \ref{app:D} and obtain that, for $\cT\gg \cL/c_o$, 
\eqref{eq:CI1} becomes independent of $\cT$ and has the form
\ba
&&\hspace*{-1.in}
\cIin (x)  \approx  \frac{ {C}_{{\bf D}} \|F\|^2}{8 \rho_o^2 c_o^2 k_o^2}
\Big[
 \sum_{j=1}^N {S}^2_{jj}  \phi_j(x)^2 \beta_{j}(\om_o)^2
 e^{\kappa_{jj}L}+
 \sum_{\stackrel{j,l=1}{j\neq l}}^N {S}_{jl}^2  \phi_l(x)^2  \beta_j(\om_o) \beta_{l}(\om_o)
e^{\kappa_{jj}L} \nonumber \\ 
&& \quad \quad +
 \sum_{j, j',l = 1}^N  {S}_{j'l}^2  \phi_l(x)^2 \beta_{j'}(\om_o) \beta_{l}(\om_o)
 \int_{-\infty}^\infty {\rm d} t \, W_{j}^{(j')}(\om_o,t,L)    \Big]. \label{eq:intContr}
\ea
This is approximately constant in $x$, at least for $\bS$ diagonally dominant, which is needed for the 
imaging function to focus at the target.

More explicitly, in the  case of an incoherent source with full aperture, where $\bS$ equals the identity matrix, 
we have in the weak  scattering regime that 
\ba
\cIin (x) \approx \frac{ {C}_{{\bf D}} \|F\|^2}{8 \rho_o^2 c_o^2 k_o^2} \Phi_2(x,x). \label{eq:intContrWeak}
\ea
Furthermore, in  the strong scattering regime, at range $L \gg L_{\rm eq}$, we obtain exactly the same expression.
Thus, the function $\cIin(x)$ is proportional to $\Phi_2(x,x)$ and is approximately constant for 
$x$ away from the boundary points $0$ and $X$. 

{
\begin{figure}[t]
\begin{center}
\includegraphics[width=2.15in]{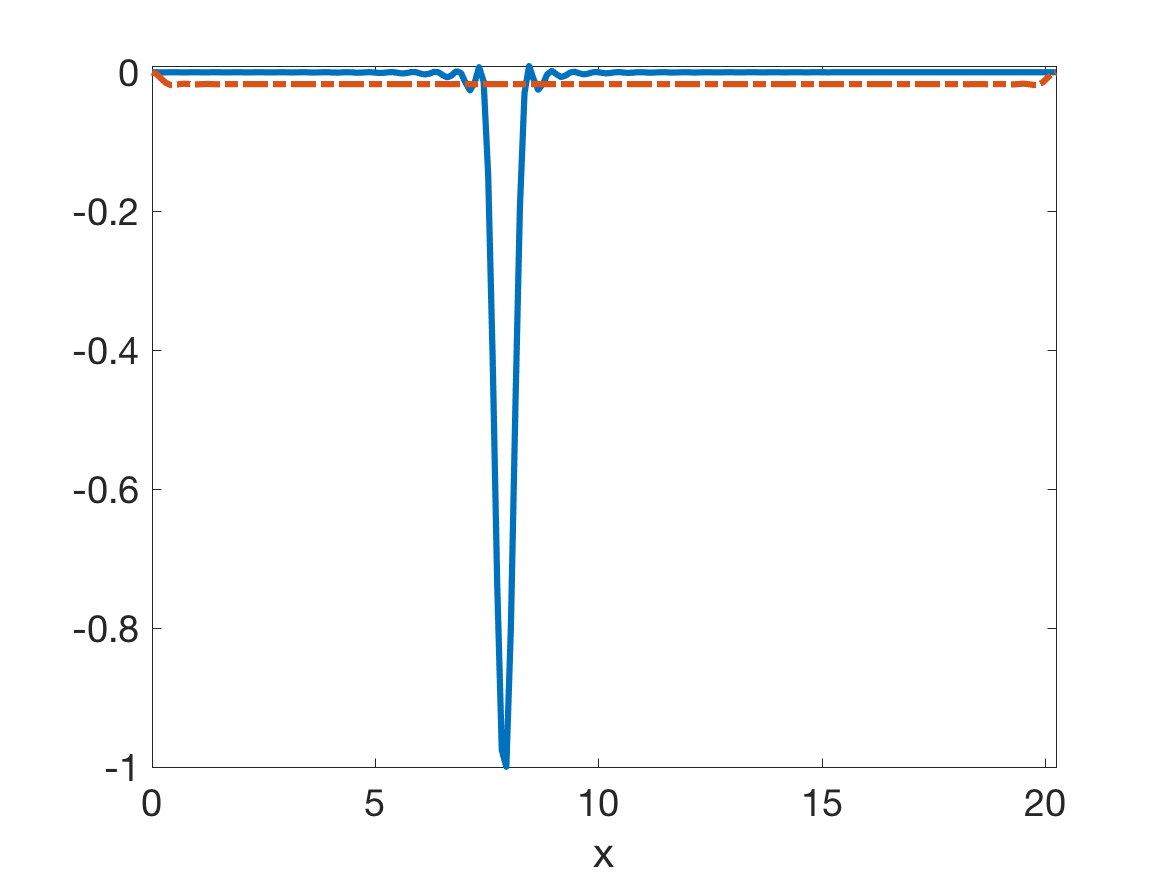}
\hspace{-0.3in}
\includegraphics[width=2.15in]{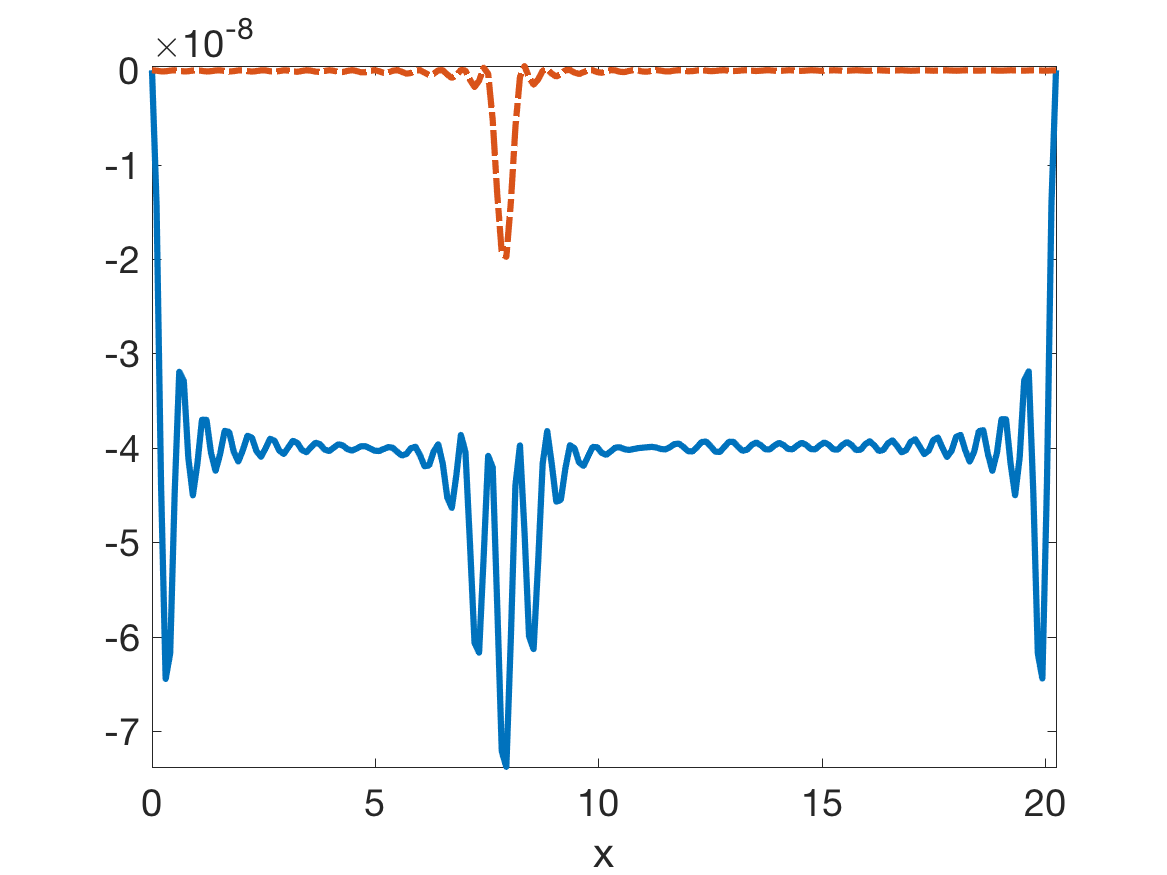}
\hspace{-0.3in}
\includegraphics[width=2.15in]{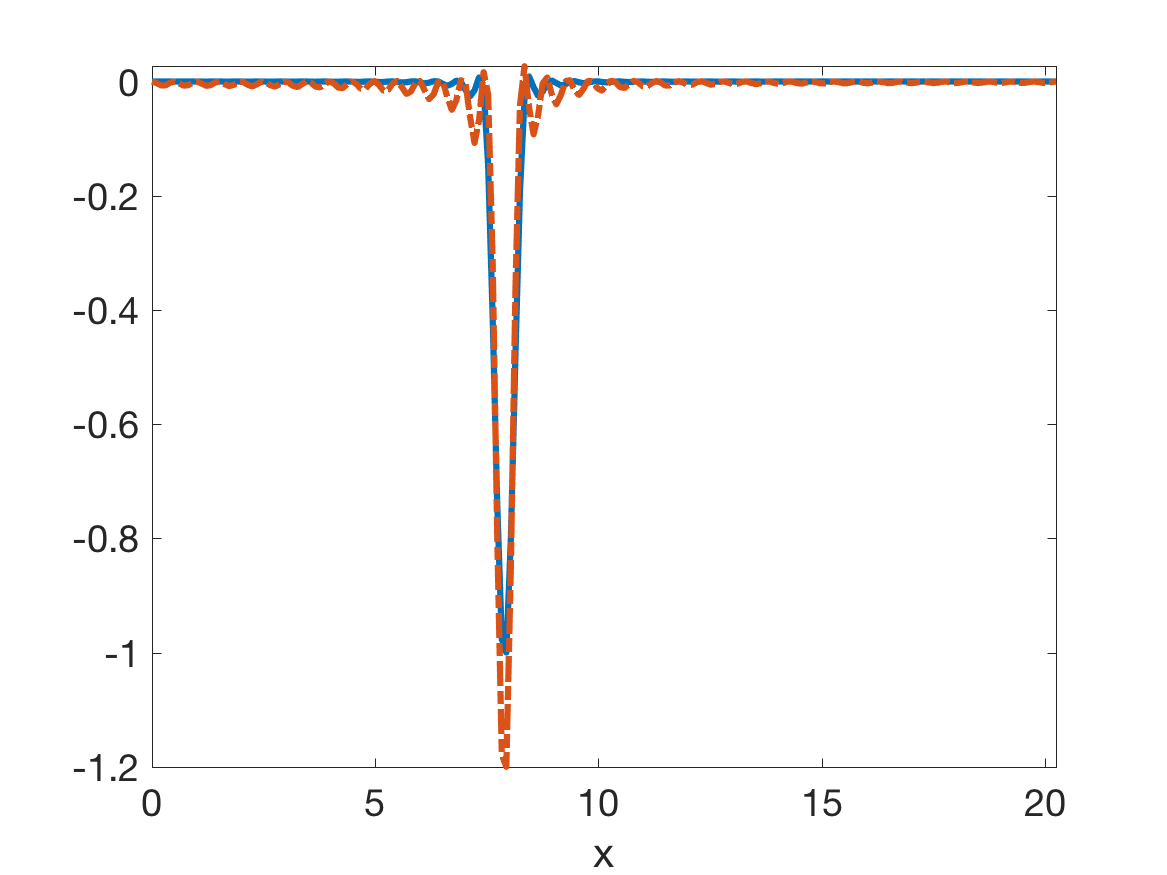}
\end{center}
\caption{
{Left: Imaging function in the homogeneous waveguide (solid blue line) and in the random waveguide
(where we cannot measure, so the reference waveguide is the  homogeneous waveguide) (red dotted line) for 
full detector and full source aperture and for $\cT > \beta'_N \cL$.  Middle: Imaging function in the homogeneous waveguide for all the recorded modes (solid line) and in the random waveguide (where we measure, so the reference waveguide is the same random waveguide)
(red dotted line) for $a_d = a_s = 0.05X$ and for $\cT > \beta'_N \cL$. Right: The same as in the middle, but for $a_d = a_s = X$.
We normalize by the absolute value at the target location 
of the image calculated in the homogeneous waveguide, with  $a_d = a_s = X$ and $\cT > \beta'_N \cL$. The abscissa 
is cross-range in units of the carrier wavelength.}}
\label{fig:Res2}
\end{figure}
}

\subsection{Imaging based on the random reference waveguide}
\label{sect:ImagR}
We present here the analysis of the imaging function \eqref{eq:Im10}, in the
waveguide with random boundary, using the random waveguide as the reference
waveguide.  
The imaging function is given by (\ref{eq:Im10}) but with 
$\cG_{11}$ replaced by 
\ba
&&\hspace*{-1.in}
\cG_{11}(\om,y,x) = \iint_{\cAs}\hspace{-0.05in} {\rm d} \xi {\rm d}
\xi' \, \theta(\xi,\xi') \hat g \big(\om,(y,L_\ep),(\xi,0)\big)
\overline{\hat{g} \big(\om,(x,L_\ep),(\xi',0)\big)},
\label{eq:defM11r}
\ea
and $\cG_{12}$ replaced by 
\ba
&&\hspace*{-1.in}
\cG_{12}(\om,y',x) = \iint_{\cAs}\hspace{-0.05in} {\rm d} \xi {\rm d}
\xi' \, \theta(\xi,\xi')\overline{ \hat g\big(\om,(y',L_\ep),(\xi,0)\big)}
(-i)\partial_z {\hat{g} \big(\om,(x,L_\ep),(\xi',0)\big)}.
\label{eq:defM12r}
\ea

If the medium is weakly scattering between the source and the target, there is no difference
compared to the case addressed in Section \ref{sect:ImagMR} and
we find that the imaging function has the form
(\ref{eq:expressimrand1}). {In this case, the quality of the image depends strongly of the source, as illustrated 
in Figure \ref{fig:Res1}.}

If the medium is strongly scattering between the source and the target, in the sense that $L$ is larger than
the equipartition distance, then the situation is 
dramatically different from the one addressed in Section \ref{sect:ImagMR}. As shown in  \ref{app:E},
we find that the imaging function has the form:
\ba
&&\hspace*{-1.in}
\cI (x) 
\approx  -\frac{k^2_o \|F\|^2 C_{{\bf D}} C_{{\bf S}}}{16 \rho_o^2c_o^2}
\iint_0^{X} {\rm d} y {\rm d}y'  {r}(y){r}(y')
  \Phi_{-1}(x,y) \Phi_{0}(x,y') \Phi_{-1}(y,y')   ,
  \label{eq:imIdeal}
\ea
where $ C_{{\bf D}}$ is given by (\ref{eq:Im13p}) and $C_{{\bf S}}$ is defined by
\ba
&&\hspace*{-1.in}
C_{{\bf S}} = \frac{\big[ \sum_{j=1}^N S_{jj} \beta_j(\om_o)\big]^2}{N(N+1)}=  \frac{1}{N(N+1)}\Big[ \iint_{\cAs} \theta(\xi,\xi') \Phi_1(\xi,\xi') {\rm d} \xi {\rm d} \xi'\Big]^2.
\label{eq:multConst}
\ea
In this case the source can be quite arbitrary. Its diameter and its coherence properties only affect the value of the constant 
$C_{{\bf S}}$. This robustness comes from the randomness of the medium that generates the suitable
illumination.  This result shows that we get an image (of the shadow) of the target, whatever the source and detector
are and that the resolution is of order $\la_o/2$.  {The only effect of the apertures is in the multiplying constants $C_{{\bf D}}$ and  $C_{{\bf S}}$ which should be large enough to make the peak observable in practice. This is illustrated in the middle and right plots in Figure \ref{fig:Res2}, for the waveguide supporting $40$ propagating modes and the point target at cross-range $0.39X$.}

\section{Summary}
\label{sect:Sum}
We analyzed a ghost imaging modality in a waveguide with random boundary, for  a penetrable scatterer (target) located 
between a partially coherent source and a detector that measures the energy flux. We considered 
a very large distance between the target and the detector, so that cumulative scattering of the waves at the random boundary
distributes the energy evenly among the propagating components (modes) of the wave. Conventional imaging 
fails in this equipartition regime, but ghost imaging is possible in three configurations:
\\
1) The source has large aperture (spans almost the entire cross-section of the waveguide) and is spatially delta-correlated, the waveguide is homogeneous,
and we correlate
the energy flux recorded at the detector with the spatially-resolved energy flux at the target distance in a reference homogeneous waveguide.
{The spatially-resolved energy flux in the reference homogeneous waveguide can be computed provided the source transmission 
is known because the waveguide modes are known.}
If the detector is large enough, imaging is possible and this is the classical ghost imaging situation.
\\
2) The source has large aperture and is delta-correlated, the waveguide is random, and we correlate
the energy flux recorded at the detector with the  spatially-resolved  energy flux at the target distance in a reference homogeneous waveguide.
Ghost imaging is possible if scattering from the source to the target is weak. 
In this configuration, scattering between the target and the detector
helps as it distributes the information over all the modes and the image becomes essentially independent of the aperture of the detector.
\\ 3) The source and the detector are arbitrary, the waveguide is random, and we correlate
the energy flux recorded at the detector with the  spatially-resolved  energy flux  at the target distance in the same random waveguide,
used as the reference waveguide.
{The spatially-resolved energy flux in the reference waveguide cannot be computed because we do not know the realization of the random medium, so it needs to be measured, or one needs to measure the transmission matrix from the source plane to the target plane
and to know the source transmission.}
Ghost imaging is then possible even if scattering from the source to the target  and from the target to the detector is strong.
This is in fact the optimal situation from the ghost imaging point of view, but it requires calibration i.e., measuring the energy flux in the waveguide in the absence  of the target.

\section*{Acknowledgments}
This material is based upon research supported in part by the Air
Force Office of Scientific Research under award FA9550-18-1-0131.

\appendix
\section{Calculation of the energy flux correlation}
\label{app:A}
The calculation of the correlation \eqref{eq:Im8} of the energy flux is
based on the expressions \eqref{eq:IW},
\eqref{eq:Born}, (\ref{eq:Im1}--\ref{eq:Im3}), the fourth-order moment formula
\ba
&&\hspace*{-1.in}
\lb {\rm d} \hat f(\om_1,x_1) {\rm d}\overline{ \hat f(\om_1',x_1')}
  {\rm d} \overline{\hat f(\om_2,x_2)} {\rm d} \hat f(\om_2',x_2') \rb
  \nonumber \\
&&\hspace*{-1.in}
 ~ ~ - \lb {\rm d} \hat f(\om_1,x_1) {\rm d}\overline{
    \hat f(\om_1',x_1')} \rb \lb {\rm d} \overline{\hat f(\om_2,x_2)}
            {\rm d} \hat f(\om_2',x_2')\rb \nonumber \\
&&\hspace*{-1.in}
~~=
            \frac{(2 \pi)^2}{B_\ep} \hat F
            \Big(\frac{\om_1}{B_\ep}\Big) \hat F
            \Big(\frac{\om_2}{B_\ep}\Big) \theta(x_1,x_2)
            \theta(x_1',x_2') \delta(\om_1-\om_2)
            \delta(\om_1'-\om_2'){\rm d} \om_1 {\rm d} \om_2{\rm d} \om_1' {\rm d} \om_2'.
\label{eq:IWf4}
\ea
derived from the Gaussian property of $\hat df(\om,x)$, and the
following lemma:
\begin{lemma}
\label{lem:1}
For a large and positive range offset $z_\ep - \Lep$, we have 
\ba
&&
\hat g\big(\om,(x,z_\ep),(\xi,0)\big) = 2
\int_0^{X} \hspace{-0.05in} {\rm d}y \, \hat
g\big(\om,(y,L_\ep),(\xi,0)\big)\hat
g\big(\om,(x,z_\ep),(y,L_\ep)\big).
\label{eq:A1}
\ea
\end{lemma}

\debproof 
Because $z_\ep - L_\ep$ is large, we can neglect the evanescent waves
and write the Green's function as a superposition of propagating modes:
\ba
\nonumber
&&\hspace*{-1.in}
\hat g\big(\om,(x,z_\ep),(y,L_\ep)\big)
 \stackrel{by~(\ref{eq:IW2})}{=}
 \sum_{j=1}^N \phi_j(x)
\frac{a_j(\om,z_\ep;(y,L_\ep))}{\sqrt{\beta_j(\om_o)}} e^{i
  \beta_j(\om+\om_o)(z_\ep - L_\ep)} \\
\nonumber
  &&\hspace*{.36in} \stackrel{by~(\ref{eq:IW8})}{=}
 \sum_{j,q=1}^N \phi_j(x) {P}_{jq}(\om,z_\ep;\Lep)
\frac{a_q(\om,L_\ep;(y,L_\ep))}{\sqrt{\beta_j(\om_o)}} e^{i
  \beta_j(\om+\om_o)(z_\ep - L_\ep)} \\
  &&\hspace*{.36in} 
  \stackrel{by~(\ref{eq:IW7})}{=}
  \frac{1}{2} 
 \sum_{j,q=1}^N \phi_j(x) \phi_q(y) {P}_{jq}(\om,z_\ep;\Lep)
\frac{\sqrt{\beta_q(\om_o)}}{\sqrt{\beta_j(\om_o)}} e^{i
  \beta_j(\om+\om_o)(z_\ep - L_\ep)} .
  \label{eq:AA1a}
\ea
We can also write by~\eqref{eq:IW2}:
\ba
\hat g\big(\om,(x,z_\ep),(\xi,0)\big) = \sum_{j=1}^N \phi_j(x)
\frac{\alpha_j(\om,z_\ep;\xi)}{\sqrt{\beta_j(\om_o)}} e^{i
  \beta_j(\om+\om_o)(z_\ep - L_\ep)},
\label{eq:AA1}
\ea
with amplitudes 
\ban
 \alpha_j(\om,z_\ep;\xi) =
a_j(\om,z_\ep;(\xi,0))e^{ i \beta_j(\om+\om_o)
  \Lep}.
\ean
These are defined by:
\ban
\alpha_j(\om,z_\ep;\xi) =\sum_{q=1}^N {P}_{jq}(\om,z_\ep;\Lep)
\alpha_q(\om,L_\ep;(\xi,0)),
\ean
in terms of the $N \times N$ propagator matrix ${\bf P}(\om,z_\ep;L_\ep)$
and the mode amplitudes of the Green's function at range $L_\ep$,
\ban
\alpha_q(\om,L_\ep;(\xi,0)) &=& a_q(\om,L_\ep;(\xi,0)) 
e^{ i \beta_q(\om+\om_o)
  \Lep}\\&=& \sqrt{\beta_q(\om_o)}
\int_0^{X} \hspace{-0.05in} {\rm d}y \, \phi_q(y) \hat
g\big(\om,(y,L_\ep),(\xi,0)\big).
\ean
This construction ensures the continuity of $\hat g$ and its range
derivative at $z_\ep = \Lep$ (see \cite{alonso2011wave} and
\cite[Section 4.4]{borcea2015inverse}).  Substituting in
\eqref{eq:AA1} we get
\ban
  &&\hspace*{-1.in} 
\hat g\big(\om,(x,z_\ep),(\xi,0)\big) =
\int_0^{X} \hspace{-0.05in} {\rm d}y \, \hat
g\big(\om,(y,L_\ep),(\xi,0)\big) \nonumber \\ 
&&\hspace*{.42in}  \times\sum_{j,q=1}^N
\sqrt{\frac{\beta_q(\om_o)}{\beta_j(\om_o)}} \phi_j(x) \phi_q(y) 
e^{i \beta_j(\om+\om_o)(z_\ep - L_\ep)} {P}_{jq}(\om,z_\ep;L_\ep),
\ean
and the result \eqref{eq:A1} follows from \eqref{eq:AA1a}.
\finproof
\begin{remark}Note that once we write the mode expansion of $\hat g$ on both sides
of equation \eqref{eq:A1}, which holds for all $x,\xi \in (0,{X})$, we obtain 
using the orthogonality of $\{\phi_j\}$ the following relation for 
the propagators
\ban
  &&\hspace*{-1.in} 
e^{i \beta_j(\om+\om_o) z_\ep} {P}_{jq}(\om,z_\ep;0) = \sum_{l=1}^N 
e^{i \beta_j(\om+\om_o) (z_\ep-\Lep)} {P}_{jl}(\om,z_\ep;\Lep) e^{i \beta_l(\om+\om_o) \Lep} 
{P}_{lq}(\om,\Lep;0).
\ean
This relation is trivial in an unperturbed waveguide, where
${P}_{jq}(\om,z_\ep;\zeta_\ep) = \delta_{jq}$.
\end{remark}

Using the expression \eqref{eq:Born} of the scattered wave,
\ba
  &&\hspace*{-1.in} 
{\rm d} \hpin(\om,x,\Lep+\cLep) \approx 2 \int_0^{X}  \hat g\big(\om,(x,\Lep + \cLep), (y,\Lep)\big) {\rm d} \hpin(\om,y,\Lep) {\rm d}y,
\label{eq:D1}
\ea
 and keeping
the terms up to $O({r}^2)$, we obtain that 
\ban
  &&\hspace*{-1.in} 
  I(t)-\Iin(t) \approx 
\frac{2}{(2 \pi)^2 \rho_o c_o}  {\rm Re}
\int_{{\cAd}} {\rm d} x'   \iint_0^{X} {\rm d} y {\rm d} y' \iint_{-\infty}^\infty 
{\rm d} \hpin(\om,y,\Lep) {\rm d}\overline{\hpin(\om',y',\Lep)} \\
  &&\hspace*{-1.in} \quad \times  e^{i(\om'-\om)t}    
\bigg\{
-2k_o 
\Big[
  \hat{G}\big(\omega,(x',\Lep+\cLep),(y,\Lep) \big) 
i\partial_z \overline{\hat{g} (\omega', (x',\Lep+\cLep),(y',\Lep)\big)} 
 \Big] r(y)\\
  &&\hspace*{-1.in}  \quad
-2k_o 
\Big[
 \hat{g}\big(\omega,(x',\Lep+\cLep),(y,L_\ep) \big) 
i\partial_z \overline{\hat{G} (\omega', (x',\Lep+\cLep),(y',\Lep)\big)} 
 \Big] r(y')
   \\
  &&\hspace*{-1.in}  \quad
  + 2k_o^3 \Big[ \int_0^{X} {\rm d} y'' \, r(y'')
    \hat{G}\big(\om,(x',\Lep+\cLep),(y'',\Lep)\big)
              \hat{G}_o\big(\om,(y'',\Lep),(y,\Lep)\big)\\
  &&\hspace*{-1.in}  \quad
  \quad \quad \quad \times
            i \partial_z \overline{\hat{g}\big(\om',(x',\Lep+\cLep),(y',\Lep)\big)} \Big] {r}(y)\\
  &&\hspace*{-1.in}  \quad
  +2k_o^3\Big[ \int_0^{X}{\rm d} y'' \, r(y'')
    \hat{g}\big(\om,(x',\Lep+\cLep),(y,\Lep)\big)
             \overline{ \hat{G}_o\big(\om',(y'',\Lep),(y',\Lep)\big)}\\
   &&\hspace*{-1.in} \quad
   \quad \quad \quad \times
            i \partial_z \overline{\hat{G}\big(\om',(x',\Lep+\cLep),(y'',\Lep)\big)} \Big] {r}(y')\\
  &&\hspace*{-1.in}  \quad
          +k_o^3 
   {r}(y)  {r}(y')  
    \hat{G}\big(\om,(x',\Lep+\cLep),(y,\Lep)\big)
   i \partial_z  \overline{\hat{G}\big(\om',(x',\Lep+\cLep),(y',\Lep)\big)}
\bigg\}.
\ean
This gives
\ban
  &&\hspace*{-1.in} 
  I(t)-\Iin(t) \approx \frac{2}{(2 \pi)^2 \rho_o c_o}  {\rm Re}   \iint_0^{X} {\rm d} y {\rm d} y'
 \iint_{-\infty}^\infty  
{\rm d} \hpin(\om,y,\Lep) {\rm d}\overline{\hpin(\om',y',\Lep)}  \\
  &&\hspace*{-1.in}  \quad\times   e^{i(\om'-\om)t}    
\bigg\{
-2k_o 
\cG_{31} (\omega,\omega',y,y') {r}(y) -2k_o \cG_{32} (\omega,\omega',y,y') {r}(y')
   \\
  &&\hspace*{-1.in}  \quad+ 2k_o^3 \Big[ \int_0^{X}{ \rm d} y'' \, r(y'')\cG_{31}\big(\omega,\omega',y'',y'\big)  
              \hat{G}_o\big(\om,(y'',\Lep),(y,\Lep)\big) \Big] {r}(y)\\
  &&\hspace*{-1.in}  \quad+2k_o^3\Big[ \int_0^{X}{\rm d} y'' \, r(y'') \cG_{32}\big(\omega,\omega',y,y''\big) 
       \overline{ \hat{G}_o\big(\om',(y'',\Lep),(y',\Lep)\big)} \Big] {r}(y')\\
  &&\hspace*{-1.in}   \quad        +k_o^3 
 \cG_2(\omega,\omega',y,y')  {r}(y)  {r}(y')\bigg\},
\ean
with $\cG_2$, $\cG_{31}$, and $\cG_{32}$ defined in (\ref{eq:defM2}--\ref{eq:defM32}).

The correlation function is obtained by substituting this expression
into \eqref{eq:Im8} and using definitions (\ref{eq:Im2}--\ref{eq:Im3})
and the moment formula \eqref{eq:IWf4}
\ba
&& \hspace*{-1.in}
\cC_T(s,x)  \approx 
\frac{2}{(2 \pi)^2 B_\ep\rho_o^2 c_o^2 } {\rm Re} \iint_{-\infty}^\infty \hspace{-0.05in} {\rm d}\om {\rm
  d}\om' \hat F\Big(\frac{\om}{B_\ep}\Big)\hat
F\Big(\frac{\om'}{B_\ep}\Big) e^{i (\omega'-\omega) s} \iint_0^{X} \hspace{-0.05in} {\rm d} y
{\rm d} y'\nonumber \\ 
&&\times  \Big[ \cG_{12}(\om',y',x)
     \cG_{11}(\om,y,x) + \overline{\cG_{12}(\om,y,x)}
     \overline{\cG_{11}(\om',y',x)} \Big]
     \nonumber\\
     &&\times 
     \Big[ 
     -2  r(y) \cG_{31}( \omega,\omega',y,y')  
     -2 r(y') \cG_{32}( \omega,\omega',y,y')  
     \nonumber \\
     &&
     \quad +2 k_o^2 r(y) \int_0^{X} {\rm d} y'' \, r(y'')\cG_{31}( \omega,\omega',y'',y')  \hat{G}_o\big(\omega,(y'',L_\ep),(y,L_\ep)\big)          \nonumber\\
     && 
      \quad     +2 k_o^2 r(y') \int_0^{X} {\rm d} y'' \, r(y'')\cG_{32}( \omega,\omega' ,y,y'') \overline{ \hat{G}_o\big(\omega',(y'',L_\ep),(y',L_\ep)\big)}
     \nonumber\\
     && \quad +   k_o^2 r(y) r(y') \cG_{2}( \omega,\omega',y,y') 
     \Big]
.
\label{eq:Im9}
\ea
Here we recall that $\hat F$ and ${r}$ are real valued and we have used the
notations (\ref{eq:defM11}--\ref{eq:defM32}). We also changed some
variables of integration.

The imaging function is given by the integral of \eqref{eq:Im9} over
$s \in (0,\tau_\ep)$, with $\tau_\ep = \cT/\ep^2$. Because
\ban
\int_0^{\tau_\ep} ds \, e^{i (\om'-\om)s} = \tau_\ep e^{i (\om'-\om)
  \frac{\tau_\ep}{2}} \mbox{sinc} \left[\frac{(\om'-\om)\tau_\ep}{2}\right],
\ean
only $O(\ep^2)$ frequency offsets contribute. Thus, we change the
frequency variables of integration
$
(\om,\om') \leadsto (\om, h),$ with $h = \frac{\om-\om'}{\ep^2},$
and approximate
$
\hat F\big(\frac{\om - \ep^2 h}{B_\ep}\big) \approx \hat
F\big(\frac{\om}{B_\ep}\big),
$
using the bandwidth scaling \eqref{eq:LR7} and assuming that $\hat F$
is continuous. The result \eqref{eq:Im10} follows.

\section{The imaging function in the unperturbed waveguide}
\label{app:B}
We begin with the
calculation of the terms
\eqref{eq:defM11}--\eqref{eq:defM32}, using the expression of the Green's functions
in the unperturbed waveguide 
\ba
\hgr\big(\om,(x,z_\ep),(\xi,\zeta_\ep)\big) &=& \hat
g\big(\om,(x,z_\ep),(\xi,\zeta_\ep)\big) 
\nonumber \\ &\approx&
\frac{1}{2} \sum_{j=1}^N \phi_j(x) \phi_j(\xi) e^{i \beta_j(\om+\om_o)
  (z_\ep-\zeta_\ep)},
\label{eq:H1} \\
\hat G\big(\om,(x,z_\ep),(\xi,\zeta_\ep)\big) &\approx &
\frac{1}{2i } \sum_{j=1}^N \frac{\phi_j(x)
  \phi_j(\xi)}{\beta_j(\om_o)} e^{i \beta_j(\om+\om_o)
  (z_\ep-\zeta_\ep)},
\label{eq:H2}
\ea
for $z_\ep> \zeta_\ep$ satisfying $z_\ep - \zeta_\ep = O(\la_o/\ep^2)$.

Substituting \eqref{eq:H1} into (\ref{eq:defM11}-\ref{eq:defM12}) we get
\ba
&&\hspace*{-1.in}
\cG_{11}(\om-\ep^2 h,y,x) \approx \frac{1}{4} \sum_{l,l'=1}^N {S}_{ll'}
\phi_{l}(x)\phi_{l'}(y) e^{- i \big\{\beta_l(\om_o) -\beta_{l'}(\om_o)+
  \om [\beta'_l(\om_o)-\beta'_{l'}(\om_o)]\big\} \frac{L}{\ep^2}}
\nonumber \\ &&\hspace*{0.5in} \times 
e^{  ih [\beta'_l(\om_o)-\beta'_{l'}(\om_o) ]L},
\label{eq:H5}\\
&& \hspace*{-1.in}
\cG_{12}(\om-\ep^2 h,y',x)  \approx \frac{1}{4} \sum_{j,j'=1}^N \beta_j(\om_o) {S}_{jj'}
\phi_j(x)\phi_{j'}(y') e^{ i \big\{\beta_j(\om_o) -\beta_{j'}(\om_o)+
  \om [\beta'_j(\om_o)-\beta'_{j'}(\om_o)]\big\} \frac{L}{\ep^2}}
\nonumber \\ &&\hspace*{0.5in} \times e^{ - ih [\beta'_j(\om_o)-\beta'_{j'}(\om_o) ]L}.
\label{eq:H6}
\ea
The substitution of (\ref{eq:H1}-\ref{eq:H2}) into \eqref{eq:defM2} gives 
\ba
&& \hspace*{-1.in}
\cG_2(\om,\om-\ep^2 h,y,y')  \approx  \frac{1}{4} \sum_{n,n'=1}^N
\frac{{D}_{nn'}}{\beta_n(\om_o)} \phi_{n}(y)
\phi_{n'}(y') e^{i \big[\beta_n(\om_o) -\beta_{n'}(\om_o)\big]
  \frac{\cL}{\ep^2}}\nonumber \\ &&\hspace*{0.5in}\times e^{ i \om
  [\beta'_n(\om_o)-\beta'_{n'}(\om_o) ] \frac{\cL}{\ep^2}+ ih
  \beta'_{n'}(\om_o)\cL},
\label{eq:H7}
\ea
and similarly, by substituting (\ref{eq:H1}-\ref{eq:H2}) into (\ref{eq:defM31} -\ref{eq:defM32})
we obtain
\ba
&& \hspace*{-1.in}
\cG_{31}(\om,\om-\ep^2 h,y,y') \approx \frac{-i}{4} \sum_{n,n'=1}^N
\frac{{D}_{nn'} \beta_{n'}(\om_o)}{\beta_{n}(\om_o)} \phi_{n}(y)
\phi_{n'}(y') e^{i [\beta_n(\om_o) -\beta_{n'}(\om_o)]
  \frac{\cL}{\ep^2}}\nonumber \\ 
  &&\hspace*{0.5in}\times e^{ i \om
  [\beta'_n(\om_o)-\beta'_{n'}(\om_o) ] \frac{\cL}{\ep^2}+ ih
  \beta'_{n'}(\om_o)\cL} ,\\
&& \hspace*{-1.in}
\cG_{32}(\om,\om-\ep^2 h,y,y') \approx \frac{i}{4} \sum_{n,n'=1}^N
 {D}_{nn'}  \phi_{n}(y)
\phi_{n'}(y') e^{i [\beta_n(\om_o) -\beta_{n'}(\om_o)]
  \frac{\cL}{\ep^2}}\nonumber \\
   &&\hspace*{0.5in}\times e^{ i \om
  [\beta'_n(\om_o)-\beta'_{n'}(\om_o) ] \frac{\cL}{\ep^2}+ ih
  \beta'_{n'}(\om_o)\cL}.
\label{eq:H8}
\ea

The imaging function \eqref{eq:Im10} involves the integral
\ba
&&
\hspace*{-1.in}
\frac{\cT}{2 \pi} \int_{-\infty}^\infty {\rm d} h \,
\mbox{sinc}\Big(\frac{h \cT}{2}\Big) e^{-i h \frac{\cT}{2} - i h [
    \beta_j'(\om_o)-\beta'_{j'}(\om_o)] L + i h \beta_{n'}'(\om_o)\cL}
\nonumber \\ 
&&
\hspace*{-1.in}= \frac{1}{2 \pi}\int_0^\cT {\rm d} t
\int_{-\infty}^\infty {\rm d} h \, e^{-i h \big\{ t -
  \beta_{n'}'(\om_o) \cL + [ \beta_j'(\om_o)-\beta'_{j'}(\om_o)] L
  \big\}} \nonumber \\ 
  \nonumber
&&
\hspace*{-1.in}=H \Big(\cT - \beta_{n'}'(\om_o) \cL + [
  \beta_j'(\om_o)-\beta'_{j'}(\om_o)] L \Big) -
  H \Big(- \beta_{n'}'(\om_o) \cL + \big[
  \beta_j'(\om_o)-\beta'_{j'}(\om_o)\big] L \Big) 
\\ 
&&
\hspace*{-1.in}=H \Big(\cT - \beta_{n'}'(\om_o) \cL + [
  \beta_j'(\om_o)-\beta'_{j'}(\om_o)] L \Big),
\label{eq:intH}
\ea 
where $H$ is the Heaviside step function and we have used the scaling
relation \eqref{eq:LR5} to conclude that 
\ban
&&
\hspace*{-1.in}
H \Big(- \beta_{n'}'(\om_o) \cL + \big[
  \beta_j'(\om_o)-\beta'_{j'}(\om_o)\big] L \Big) = 0, \quad \forall
j,j',n' = 1, \ldots, N.
\ean
We also have the integral
\ban
&&
\hspace*{-1.in}\frac{1}{2 \pi B_\ep}\int_{-\infty}^\infty {\rm d}
  \om \, \hat F^2\Big(\frac{\om}{B_\ep}\Big) e^{-i \frac{\om}{\ep^2}
    \big\{ [\beta'_{n'}(\om_o)-\beta'_n(\om_o)] \cL - [
      \beta_j'(\om_o)-\beta'_{j'}(\om_o) +
      \beta'_{l'}(\om_o)-\beta'_l(\om_o)] L\big\}} \nonumber \\ 
&&
\hspace*{-1.in}=
  F\star F \left( \frac{B_\ep
    \cL}{\ep^2}\Big[\beta'_{n'}(\om_o)-\beta'_n(\om_o)- \big(
    \beta_j'(\om_o)-\beta'_{j'}(\om_o) +
    \beta'_{l'}(\om_o)-\beta'_l(\om_o)\big) \frac{L}{\cL} \Big]
  \right),
%\label{eq:intom}
\ean
where $\star$ denotes convolution. 

Gathering the results (\ref{eq:H5}--\ref{eq:intH}) and substituting
into \eqref{eq:Im10} we obtain that the imaging function is the sum of
three terms:
\ba
\cI_{\frac{\cT}{\ep^2}}(x)\approx \sum_{j=1}^3 \cI_{\frac{\cT}{\ep^2},j}(x).
\ea
The first term is
\ba
\nonumber
&&
\hspace*{-1.in}
\cI_{\frac{\cT}{\ep^2},1}(x) = \frac{1}{16
\rho_o^2c_o^2 } {\rm Re}  \, i \iint_0^{X} {\rm d} y {\rm d} y'
 \sum_{j,j',l,l',n,n'=1}^N  (\beta_j +\beta_l)
  {S}_{jj'} \phi_j(x) \phi_{j'}(y') 
 \\
&&
\hspace*{-0.2in}
\times
{S}_{ll'} \phi_l(x) \phi_{l'}(y)e^{i
  [\beta_j(\om_o)-\beta_{j'}(\om_o) -\beta_l(\om_o)+\beta_{l'}(\om_o)] \frac{L}{\ep^2}} 
  \nonumber
 \\ 
 &&
\hspace*{-0.2in} \times 
 {D}_{nn'} \phi_n(y) \phi_{n'}(y') \left[ \frac{\beta_{n'}(\om_o)}{\beta_n(\om_o)} {r}(y) 
 - {r}(y')\right]
 e^{i  [\beta_n(\om_o)-\beta_{n'}(\om_o)] \frac{\cL}{\ep^2}}
\nonumber \\ 
&&
\hspace*{-0.2in}
\times F\star F \left( \frac{B_\ep
  \cL}{\ep^2}\Big[\beta'_{n'}(\om_o)-\beta'_n(\om_o)- \big(
  \beta_j'(\om_o)-\beta'_{j'}(\om_o) +
  \beta'_{l'}(\om_o)-\beta'_l(\om_o)\big) \frac{L}{\cL} \Big]
\right)\nonumber \\
&&
\hspace*{-0.2in}
\times H \Big(\cT - \beta_{n'}'(\om_o) \cL + \big[
  \beta_j'(\om_o)-\beta'_{j'}(\om_o)\big] L \Big).
\label{eq:ImHom1}
\ea
The second term is 
\ba
\nonumber
&&
\hspace*{-1.in}
\cI_{\frac{\cT}{\ep^2},2}(x) 
=
 \frac{k_0^2}{16
\rho_o^2c_o^2 } {\rm Re} \, i \iint_0^{X} {\rm d} y {\rm d} y'
 \sum_{j,j',l,l',n,n'=1}^N  (\beta_j +\beta_l)
  {S}_{jj'} \phi_j(x) \phi_{j'}(y') 
 \\
&&\hspace*{-0.2in}\times
{S}_{ll'} \phi_l(x) \phi_{l'}(y)e^{i
  [\beta_j(\om_o)-\beta_{j'}(\om_o) -\beta_l(\om_o)+\beta_{l'}(\om_o)] \frac{L}{\ep^2}} 
  \nonumber
 \\
   \nonumber 
&&\hspace*{-0.2in} \times 
 {D}_{nn'} 
 \Big[
 - \frac{\beta_{n'}(\om_o)}{\beta_n(\om_o)} \int_0^{X} \phi_n(y'') \phi_{n'}(y') \hat{G}_o\big(\omega,(y'',\Lep),(y,\Lep)\big)
  {r}(y'') {\rm d} y''  {r}(y)  \\
   \nonumber 
&&\hspace*{-0.2in} 
\quad +\int_0^{X} \phi_n(y) \phi_{n'}(y'') \overline{\hat{G}_o\big(\omega,(y'',\Lep),(y',\Lep)\big)}
  {r}(y'') {\rm d} y''  {r}(y')  
  \Big]
 e^{i  [\beta_n(\om_o)-\beta_{n'}(\om_o)] \frac{\cL}{\ep^2}}
\nonumber \\ 
&&\hspace*{-0.2in}\times F\star F \left( \frac{B_\ep
  \cL}{\ep^2}\Big[\beta'_{n'}(\om_o)-\beta'_n(\om_o)- \big(
  \beta_j'(\om_o)-\beta'_{j'}(\om_o) +
  \beta'_{l'}(\om_o)-\beta'_l(\om_o)\big) \frac{L}{\cL} \Big]
\right)\nonumber \\
&&\hspace*{-0.2in}\times H \Big(\cT - \beta_{n'}'(\om_o) \cL + \big[
  \beta_j'(\om_o)-\beta'_{j'}(\om_o)\big] L \Big).
\label{eq:ImHom2}
\ea
with $\hat G_o\big(\om,(y'',\Lep),(y',\Lep)\big)$ given in
\eqref{eq:Born1}. 
The third term is
\ba
\nonumber
&&
\hspace*{-1.in}
\cI_{\frac{\cT}{\ep^2},3}(x)
= \frac{k_o^2}{32
\rho_o^2c_o^2 } {\rm Re} \iint_0^{X} {\rm d} y {\rm d} y'
 \sum_{j,j',l,l',n,n'=1}^N  (\beta_j +\beta_l)
  {S}_{jj'} \phi_j(x) \phi_{j'}(y') 
 \\
&&\hspace*{-0.2in}\times
{S}_{ll'} \phi_l(x) \phi_{l'}(y)e^{i
  [\beta_j(\om_o)-\beta_{j'}(\om_o) -\beta_l(\om_o)+\beta_{l'}(\om_o)] \frac{L}{\ep^2}} 
  \nonumber
 \\ 
&&\hspace*{-0.2in} \times 
 {D}_{nn'} \phi_n(y) \phi_{n'}(y')  \frac{1}{\beta_n(\om_o)} {r}(y)  {r}(y') 
 e^{i  [\beta_n(\om_o)-\beta_{n'}(\om_o)] \frac{\cL}{\ep^2}}
\nonumber \\ 
&&\hspace*{-0.2in}\times F\star F \left( \frac{B_\ep
  \cL}{\ep^2}\Big[\beta'_{n'}(\om_o)-\beta'_n(\om_o)- \big(
  \beta_j'(\om_o)-\beta'_{j'}(\om_o) +
  \beta'_{l'}(\om_o)-\beta'_l(\om_o)\big) \frac{L}{\cL} \Big]
\right)\nonumber \\
&&\hspace*{-0.2in}\times H \Big(\cT - \beta_{n'}'(\om_o) \cL + \big[
  \beta_j'(\om_o)-\beta'_{j'}(\om_o)\big] L \Big).
\label{eq:ImHom3}
\ea

In all these expressions we note that since $F$ has finite support and
$B_\ep/\ep^2 \to \infty $ as $\ep \to 0$, only the terms that give
\ban
\beta'_{n'}(\om_o)-\beta'_n(\om_o)- \big(
  \beta_j'(\om_o)-\beta'_{j'}(\om_o) +
  \beta'_{l'}(\om_o)-\beta'_l(\om_o)\big) \frac{L}{\cL} = 0
\ean
contribute. Moreover, since $\cL \gg L$, we have two cases: (1) For $n = n'$ and $j=j'$ and $l = l'$.
 (2) For $n = n'$ and $j=l$ and $j' = l'$.
After
considering these two cases in \eqref{eq:ImHom1}--\eqref{eq:ImHom3},
and using that
\ban
F\star F(0) = \int_{-\infty}^\infty \frac{{\rm d} w}{2 \pi} \hat F(w)^2 = \|F\|^2,
\ean
and 
\ban
\mbox{Re} \Big[ i \hat
    G_o\big(\om,(y'',\Lep),(y',\Lep)\big) \Big] \approx  \frac{1}{2}\sum_{q=1}^N
\frac{\beta_q(y')\beta_q(y'')}{\beta_q(\om_o)} = \frac{1}{2}\Phi_{-1}(y',y'') ,
\ean
we get
\ba
\nonumber
&&\hspace*{-1.in}
\cI_{\frac{\cT}{\ep^2}}(x) 
= \frac{k_o^2 \|F\|^2}{32
\rho_o^2c_o^2 }   \iint_0^{X} {\rm d} y {\rm d} y'
\sum_{j,l=1}^N  \big(\beta_j(\om_o) +\beta_l(\om_o)\big)
  {S}_{jj} \phi_j(x) \phi_{j}(y') 
{S}_{ll} \phi_l(x) \phi_{l}(y)
  \nonumber
 \\
\nonumber
&& \hspace*{-0.2in}\times 
\sum_{n=1}^N H\big(\cT-\beta_n'(\om_o)\cL\big) D_{nn}\Big[  \frac{\phi_n(y) \phi_{n}(y') }{\beta_n(\om_o)}   {r}(y)  r(y')
\nonumber \\
&&\hspace*{-0.2in}\hspace{1.3in}-2 r(y')\int_0^{X} {\rm d} y'' \, r(y'') \phi_n(y)\phi_n(y'')\Phi_{-1}(y'',y')  \Big]
\nonumber\\
&& \hspace*{-0.2in}+
\frac{k_o^2}{16\rho_o^2c_o^2 }   \iint_0^{X} {\rm d} y {\rm d} y'
\sum_{{j'\ne j}}^N
  \beta_j(\om_o)
  {S}_{jj'}^2 \phi_j^2(x) \phi_{j'}(y)  \phi_{j'}(y') 
  \nonumber
 \\
 &&\hspace*{-0.2in} \times  \sum_{n=1}^N
H\big(\cT-\beta_n'(\om_o)\cL + (\beta'_j(\om_o)-\beta_{j'}'(\om_o))L\big)D_{nn}\Big[  \frac{\phi_n(y) \phi_{n}(y') }{\beta_n(\om_o)}  {r}(y) 
r(y')
\nonumber \\
&&\hspace*{-0.2in}\hspace{1.1in}-2 r(y') \int_0^{X}{\rm d} y'' \, r(y'')  \phi_n(y)\phi_n(y'')\Phi_{-1}(y',y'')  \Big].
\label{eq:ImHom}
\ea
The focusing of \eqref{eq:ImHom} is mostly dependent on the source aperture and coherence, which define the matrix $\bS$.
The time parameter $\cT$ does not play a significant role, as long as its greater than  $\beta_1'(\om_o)
\cL$, the scaled travel time of the fastest propagating mode between the target and detector. This is required 
to have at least one term in the sum over $n$ and thus obtain a non-trivial image.  
The imaging function  simplifies when $\cT \gg \cL/c_o$, or more precisely $\cT > \beta'_N(\om_o) \cL$, because all the Heaviside functions 
in \eqref{eq:ImHom} equal 1. The expression \eqref{eq:Im11} is obtained from \eqref{eq:ImHom} at such large $\cT$.

\section{The imaging function in the random waveguide}
\label{app:C}
Here we derive the expression of the imaging function in the random
waveguide, in the case of the reference waveguide with unperturbed
boundary. 

The expectations of $\cG_2$,  $\cG_{31}$, and $\cG_{32}$ are obtained from definitions 
\eqref{eq:defM2}--\eqref{eq:defM32} and the expressions
\ban
&& \hspace*{-1.in}
\hat G\big(\om,(x',\Lep+\cLep),(y,\Lep)\big) \approx \frac{1}{2i}
\sum_{j,q=1}^N
\frac{\phi_j(x')\phi_q(y)}{\sqrt{\beta_j(\om_o)\beta_q(\om_o)}} e^{i
  \beta_j(\om+\om_o) \cLep}
{P}_{jq}(\om,\Lep+\cLep;\Lep) , \\
&& \hspace*{-1.in}
\hat g\big(\om,(x',\Lep+\cLep),(y,\Lep)\big) 
 \approx  \frac{1}{2}
\sum_{j,q=1}^N
\sqrt{\frac{\beta_q(\om_o)}{\beta_j(\om_o)}}\phi_j(x')\phi_q(y) e^{i
  \beta_j(\om+\om_o) \cLep} 
{P}_{jq}(\om,\Lep+\cLep;\Lep),
\ean
of the Green's functions evaluated at $z_\ep = \Lep+\cLep$ and $\zeta_\ep = \Lep$. We have 
\ba
&& \hspace*{-1.in}
\EE\Big[ \cG_2(\om,\om ,y,y') \Big]  \approx 
 \frac{1}{4}
\sum_{n,m=1}^N \frac{\phi_n(y) \phi_n(y')}{\beta_n(\om_o)} 
{D}_{mm}   \EE[ |{P}_{mn}(\omega,\Lep+\cLep;\Lep)|^2]
\nonumber \\
&& \hspace*{0.4in} + \mbox{coherent part},
\label{eq:meanG2}
\ea
\ba
&& \hspace*{-1.in}
\EE\Big[ \cG_{31}(\om,\om ,y'',y') \Big]  \approx \frac{1}{4i}
\sum_{n,m=1}^N \phi_n(y'') \phi_n(y') 
 {D}_{mm} 
   \EE[ |{P}_{mn}(\omega,\Lep+\cLep;\Lep)|^2]
 \nonumber \\
&&\hspace*{0.4in} + \mbox{coherent part},
\label{eq:meanG31}
\ea
and
\ba
&& \hspace*{-1.in}
\EE\Big[ \cG_{32}(\om,\om,y,y'') \Big] \approx - \frac{1}{4i}
\sum_{n,m=1}^N \phi_n(y'') \phi_n(y) 
 {D}_{mm} 
   \EE[ |{P}_{mn}(\omega,\Lep+\cLep;\Lep)|^2] \nonumber \\
&&\hspace*{0.4in} + \mbox{coherent part}.
\label{eq:meanG32}
\ea

The coherent parts in these expressions are $O(e^{\kappa_{11} \cL})$ because as shown
in \cite{alonso2011wave},
\ban
\mbox{Re}[\kappa_{jj'}] \le  \kappa_{11} < 0, \quad \forall j,
j' = 1, \ldots, N.
\ean
We are interested in a very large range offset $
\cL \gg L_{\rm eq} $, 
where $L_{\rm eq}$ is the equipartition distance defined in
\cite[Section 20.3.3]{fouque07}. At such range 
scattering at the random boundary distributes the energy evenly among
the propagating modes so that we get
\ba
  \EE[ |{P}_{mn}(\omega,\Lep+\cLep;\Lep)|^2] \approx \frac{1}{N}
\ea
and consequently
\ba
\EE\Big[ \cG_2(\om,\om ,y,y') \Big] &\approx& \frac{C_{\bD}}{4}  \Phi_{-1}(y,y')
\label{eq:meanaG2}
 ,\\
\EE\Big[ \cG_{31}(\om,\om ,y'',y') \Big] &\approx&
 \frac{C_{\bD}}{4i} \Phi_{0}(y'',y'),
\label{eq:meanaG31}
 \\ 
\EE\Big[ \cG_{32}(\om,\om,y,y'') \Big] &\approx& -  \frac{C_{\bD}}{4i}  \Phi_{0}(y'',y),
\label{eq:meanaG32}
\ea
with $C_{\bD}$ defined in \eqref{eq:Im13p}.

The expectation of $\cG_{11} \cG_{12}$ 
is obtained from definitions \eqref{eq:defM11}-\eqref{eq:defM12} and
\ban
\hat g\big(\om,(y,\Lep) ,(\xi,0)\big) & \approx& \frac{1}{2}
\sum_{j,l=1}^N
\sqrt{\frac{\beta_l(\om_o)}{\beta_j(\om_o)}}\phi_j(y)\phi_l(\xi) e^{i
  \beta_j(\om+\om_o) \Lep} 
{P}_{jq}(\om,\Lep ; 0), \\
\hat g^{({\rm r})}\big(\om,(y,\Lep) ,(\xi,0)\big) & \approx& \frac{1}{2}
\sum_{j,l=1}^N
\sqrt{\frac{\beta_l(\om_o)}{\beta_j(\om_o)}}\phi_j(y)\phi_l(\xi) e^{i
  \beta_j(\om+\om_o) \Lep} ,
\ean
which give
\ba
\nonumber
&& \hspace*{-1.in}\EE\big[ \cG_{11} (\omega,y,x)\cG_{12}(\omega,y',x) \big]
\\
\nonumber
&& \hspace*{-1.in}\approx
\frac{1}{16} \sum_{j,j',l,l'=1}^N {S}_{jl} \phi_j(y) \phi_l(x) {S}_{j'l'}\phi_{j'}(y')\phi_{l'}(x)\beta_{l'}(\om_o)
 \EE[ {P}_{jj}(\om,\Lep;0)  \overline{{P}_{j'j'}(\om,\Lep;0)}] \\
\nonumber
&& \hspace*{-1.in}\quad \quad\times
 e^{i [\beta_j-\beta_l-\beta_{j'}+\beta_{l'}](\om_o) \frac{L}{\ep^2}+ i \omega [\beta_j'-\beta_l'-\beta_{j'}'+\beta_{l'}'](\om_o)\frac{L}{\ep^2}}
 \\
&& \hspace*{-1.in} \quad +
\frac{1}{16} \sum_{{j ,j',l,l' = 1}}^N  {S}_{j'l}{S}_{j'l'}  \phi_l(x)\phi_{l'}(x) \phi_j(y) \phi_j(y') \frac{\beta_{j'}(\om_o) \beta_{l'}(\om_o)}{\beta_j(\om_o)}
 \EE [ |{P}_{jj'}(\om,\Lep;0)|^2 ]  
 \nonumber\\
&& \hspace*{-1.in}\quad\quad \times
  e^{i [-\beta_l +\beta_{l'}](\om_o) \frac{L}{\ep^2}+ i \omega [ -\beta_l'+\beta_{l'}'](\om_o)\frac{L}{\ep^2}}.
\ea
Recalling the moment formula \eqref{eq:MomProp}, we see that the first term in the right-hand side is the coherent contribution, with
\ba
\EE[ {P}_{jj}(\om,\Lep;0)  \overline{{P}_{j'j'}(\om,\Lep;0)}] = e^{\kappa_{jj'} L} ,
\ea
while the second term is the incoherent contribution, with 
\ba
\EE[ |{P}_{jj'}(\om,\Lep;0)|^2 ]  = \int_{-\infty}^\infty {\rm d} t \, W_{j}^{(j')}(\om_o,t,L) .
\ea

As in the calculation in  \ref{app:B}, for the homogeneous waveguide, the integral in $\omega$ selects the terms
$(j=l,j'=l')$ and $(j=j',l=l')$ for the coherent contribution and 
 $(l=l')$ for the incoherent contribution.  More precisely, we have 
 \ba
\nonumber
&& \hspace*{-1.in}
\frac{1}{2\pi B_\ep} \int_{-\infty}^\infty \hspace{-0.05in} {\rm d}\om \, \hat
F\Big(\frac{\om}{B_\ep}\Big)^2
 \EE\big[ \cG_{11} (\omega,y,x)\cG_{12}(\omega,y',x) \big]
\\
\nonumber
&& \hspace*{-1.in}
\approx
\frac{\|F\|^2}{16} \Big\{\sum_{j,j'=1}^N {S}_{jj} \phi_j(y) \phi_j(x) {S}_{j'j'}\phi_{j'}(y')\phi_{j'}(x)\beta_{j'}(\om_o)
 e^{\kappa_{jj'}L} \\ 
\nonumber
&& \hspace*{-0.4in}\ +
\sum_{\stackrel{j,l=1}{j\neq l}}^N {S}_{jl}^2 \phi_j(y) \phi_l(x)^2  \phi_{j}(y') \beta_{l}(\om_o)
 e^{\kappa_{jj}L} \\ 
&& \hspace*{-0.4in} +
 \sum_{{j,j',l=1}}^N  {S}_{j'l}^2  \phi_l(x)^2 \phi_j(y) \phi_j(y') \frac{\beta_{j'}(\om_o) \beta_{l}(\om_o)}{\beta_j(\om_o)}
  \int_{-\infty}^\infty {\rm d} t \, W_{j}^{(j')}(\om_o,t,L) \Big\} .
 \label{eq:appC10}
\ea
Substituting these expressions and (\ref{eq:meanaG2}-\ref{eq:meanaG32}) into (\ref{eq:meancI1}), we obtain after some 
straightforward calculations that 
\ba
&& \hspace*{-1.in}
\cI(x) \approx  \frac{k_o^2 C_{\bD} \|F\|^2}{32 \rho_o^2 c_o^2} \iint_0^X {\rm d} y {\rm d} y' 
\Big[\sum_{j,j'=1}^N S_{jj} \phi_j(x) \phi_j(y) S_{j'j'} \beta_{j'}(\om_o)\phi_{j'}(x) \phi_{j'}(y')e^{\kappa_{jj'}L} 
\nonumber \\ 
&& \hspace*{-.3in}+ \sum_{j,j'=1}^N S_{jj} \phi_j(x) \phi_j(y') S_{j'j'} \beta_{j'}(\om_o)\phi_{j'}(x) \phi_{j'}(y)e^{\overline{\kappa_{jj'}}L} 
\nonumber \\
&& \hspace*{-.3in}+ 2 \sum_{\stackrel{j,l = 1}{j \ne l}}^N S^2_{jl} \beta_l(\om_o) \phi_l(x)^2 \phi_j(x) \phi_j(y') e^{\kappa_{jj}L} \nonumber \\
&& \hspace*{-.3in}+ 2 \sum_{{j,j',l = 1}}^N S^2_{j'l} \frac{\beta_l(\om_o) \beta_{j'}(\om_o)}{\beta_j(\om_o)} 
\phi_l^2(x) \phi_j(y) \phi_j(y')  \int_{-\infty}^\infty {\rm d} t \, W_{j}^{(j')}(\om_o,t,L) \Big] \nonumber \\
&& \hspace*{-.3in}\times
\Big[ 
-2 r(y) \int_0^X {\rm d} y'' r(y'') \Phi_0(y',y'') \Phi_{-1}(y,y'') 
 + r(y) r(y') \Phi_{-1}(y,y')\Big].
 \label{eq:ImagRand}
\ea

When scattering is weak between $0$ and $L$, we have $|\kappa_{jj'}|L \ll 1$ and the incoherent term, 
proportional to the integral of $W_j^{(j')}(\om_o,t,L)$, is negligible. The expression \eqref{eq:expressimrand1} of 
the imaging function is obtained from \eqref{eq:ImagRand}.

When scattering is strong between $0$ and $L$ i.e., $|\kappa_{jj'}|L > 1$ for all $j,j' = 1, \ldots, N$, the coherent terms in 
 \eqref{eq:ImagRand} are negligible, and the imaging function reduces to \eqref{eq:expressimrand2}.

\section{Derivation of the contribution of the incident energy flux} 
\label{app:D}
The incident energy flux is calculated using definition \eqref{eq:IW9} and the approximation (\ref{eq:D1}).
We obtain that 
\ban
&& \hspace*{-1.in}
\Iin(t) \approx 
\frac{8}{(2 \pi)^2 \rho_o c_o k_o}  {\rm Re}
\int_{{\cAd}} {\rm d} x'    \iint_0^{X} {\rm d} y {\rm d} y'
\iint_{-\infty}^\infty    {\rm d} \hpin(\om,y,\Lep) {\rm d}\overline{\hpin(\om',y',\Lep)} \\
&&\times  e^{i(\om'-\om)t}  
 \hat{g}\big(\omega,(x',\Lep+\cLep),(y,\Lep) \big) 
i\partial_z \overline{\hat{g} (\omega', (x',\Lep+\cLep),(y',\Lep)\big)} 
,
\ean
and proceeding as in  \ref{app:A} we get the following expression of \eqref{eq:CI1} 
\ba
&& \hspace*{-1.in}
 \cIin_{\frac{\cT}{\ep^2}}(x) \approx   \frac{8 \cT}{(2 \pi)^2 \rho_o^2 c_o^2 k_o^2
  B_\ep}
  {\rm Re} \int_{-\infty}^\infty \hspace{-0.05in} {\rm d}\om \, \hat{F}^2\Big(\frac{\om}{B_\ep}\Big)  \int_{-\infty}^\infty {\rm d}h \,
\mbox{sinc} \left(\frac{h \cT}{2}\right)e^{-i h \frac{\cT}{2}} 
\iint_0^{X} {\rm d} y 
{\rm d}y'\nonumber \\ 
&&\times 
     \Big[ \cG_{12}(\om-\ep^2 h,y',x)
     \cG_{11}(\om,y,x) + \overline{\cG_{12}(\om,y,x)}
     \overline{\cG_{11}(\om-\ep^2 h,y',x)} \Big]
     \nonumber  \\
&&\times 
    \cG_4(\om,\om-\ep^2 h,y,y'),
    \label{eq:D2}
\ea
with 
\ba
&& \hspace*{-1.in}
\cG_4(\om,\om',y,y') = \int_{\cAd} \hspace{-0.06in}{\rm d} x' \, \hat
g\big(\om,(x',\Lep+\cLep),(y,\Lep)\big) i \partial_z \overline{\hat
  g\big(\om',(x',\Lep+\cLep),(y',\Lep)\big)}.
\label{eq:D3}
\ea
The same argument as in section \ref{sect:ImagMR}, shows that that \eqref{eq:D2} is self-averaging, 
and it has the expression
\ba
&& \hspace*{-1.in}
\cIin_{\frac{\cT}{\ep^2}}(x)  \approx   \frac{8 \cT}{(2 \pi)^2 \rho_o^2 c_o^2 k_o^2
  B_\ep} {\rm Re} 
 \int_{-\infty}^{\infty} {\rm d} \om\,  \hat F^2\Big(\frac{\om}{B_\ep}\Big) 
\int_{-\infty}^\infty {\rm d} h \, \mbox{sinc} \Big(\frac{h \cT}{2}\Big) e^{-i h \frac{\cT}{2}} 
\iint_0^{X} {\rm d} y 
{\rm d}y' \nonumber \\
&&  \hspace*{-.3in}\times   \EE \Big[ \cG_{12}(\om-\ep^2 h,y',x)
     \cG_{11}(\om,y,x) + \overline{\cG_{12}(\om,y,x)}
     \overline{\cG_{11}(\om-\ep^2 h,y',x)} \Big] \nonumber\\
&& \hspace*{-.3in}\times 
 \EE\Big[ \cG_4(\om,\om-\ep^2 h,y,y')\Big],
\label{eq:D2p}
\ea
and for $\cT\gg \cL/c_o$, it becomes independent of $\cT$, 
\ba
&& \hspace*{-1.in}
\cIin (x) \approx \frac{8 }{2 \pi  \rho_o^2 c_o^2 k_o^2
  B_\ep}{\rm Re} 
 \int_{-\infty}^{\infty} {\rm d} \om\,  \hat F^2\Big(\frac{\om}{B_\ep}\Big) 
\iint_0^{X} {\rm d} y 
{\rm d}y' \nonumber \\
&& \hspace*{-.3in}\times   \EE \Big[ \cG_{12}(\om,y',x)
     \cG_{11}(\om,y,x) + \overline{\cG_{12}(\om,y,x)}
     \overline{\cG_{11}(\om,y',x)} \Big] 
 \EE\Big[ \cG_4(\om,\om,y,y')\Big].
\label{eq:D2pp}
\ea
The first square bracket was computed in  \ref{app:C} and 
\ba
&& \hspace*{-1.in}
\EE\Big[ \cG_4(\om,\om ,y,y') \Big] \approx 
 \frac{1}{4}
\sum_{n,m=1}^N  \phi_n(y) \phi_n(y') \beta_n(\om_o) 
{D}_{mm}   \EE[ |{P}_{mn}(\omega,\Lep+\cLep;\Lep)|^2] 
\nonumber \\
&&\hspace*{.4in}
+ \mbox{coherent part}.
\ea
When $\cL$ is much larger than the equipartition distance,
this becomes
\ba
&& \hspace*{-1.in}
\EE\Big[ \cG_4(\om,\om ,y,y') \Big] \approx
 \frac{1}{4N} \sum_{m=1}^N 
{D}_{mm}  
\sum_{n=1}^N \beta_n(\om_o)  \phi_n(y) \phi_n(y')   = \frac{C_{\bD}}{4}  \Phi_1(y,y')
\ea
with ${C}_{{\bf D}}$ defined in \eqref{eq:Im13p}.
This gives with (\ref{eq:appC10}):
\ba
&& \hspace*{-1.in}\cIin (x) \approx \frac{ {C}_{{\bf D}} \|F\|^2}{8 \rho_o^2 c_o^2 k_o^2}
\iint_0^{X} {\rm d} y 
{\rm d}y'  \Phi_1(y,y') \nonumber \\
\nonumber
&& \hspace*{-.2in}\times  
\Big[
 \sum_{j,j'=1}^N {S}_{jj} \phi_j(y) \phi_j(x) {S}_{j'j'}\phi_{j'}(y')\phi_{j'}(x)\beta_{j'}(\om_o)
 e^{\kappa_{jj'}L}\\ 
\nonumber
&& \hspace*{-.2in}\quad +
 \sum_{\stackrel{j,l=1}{j\neq l}}^N {S}_{jl}^2  \phi_l(x)^2  \phi_j(y)\phi_{j}(y') \beta_{l}(\om_o)
e^{\kappa_{jj}L} \\ 
&& \hspace*{-.2in} \quad +
 \sum_{j, j',l = 1}^N  {S}_{j'l}^2  \phi_l(x)^2 \phi_j(y) \phi_j(y') \frac{\beta_{j'}(\om_o) \beta_{l}(\om_o)}{\beta_j(\om_o)}
 \int_{-\infty}^\infty {\rm d} t \, W_{j}^{(j')}(\om_o,t,L)    \Big],
\ea
and using that 
\ban
\iint_0^{X} {\rm d} y 
{\rm d}y'  \Phi_1(y,y')  \phi_j(y) \phi_{j'}(y') = \delta_{jj'} \beta_j(\om_o),
\ean
we obtain equation \eqref{eq:intContr}.

\section{The imaging function in the random waveguide with random reference waveguide}
\label{app:E}
Here we derive the expression of the imaging function in the random
waveguide, in the case where the reference waveguide is the same random waveguide.
The only difference with respect to  \ref{app:C} is that we need to revisit the calculations
of the moments of the form $\EE[\cG_{11}\cG_{12}]$ with the new expressions
(\ref{eq:defM11r}-\ref{eq:defM12r}) for $\cG_{11}$ and $\cG_{12}$.

We have
\ba
\nonumber
&& \hspace*{-1.in}\EE\big[\cG_{11}(\omega,y,x)\cG_{12}(\omega,y',x)\big]
= \frac{1}{16} \sum_{j_1,l_1,j_1',l_1',j_2,l_2,j_2',l_2'=1}^N
\iint_{\cAs}  {\rm d} \xi_1{\rm d} \xi_1'{\rm d} \xi_2{\rm d} \xi_2'\theta(\xi_1,\xi_1')\theta(\xi_2,\xi_2') \\
\nonumber
&& \hspace*{-.8in}\times
\phi_{j_1}(y) \phi_{l_1}(\xi_1) \phi_{j_1'}(x) \phi_{l_1'}(\xi_1') 
\phi_{j_2}(y') \phi_{l_2}(\xi_2) \phi_{j_2'}(x) \phi_{l_2'}(\xi_1') 
 \frac{\sqrt{\beta_{l_1}\beta_{l_1'} \beta_{l_2}\beta_{l_2'}(\om_o)}}
 {\sqrt{\beta_{j_1}\beta_{j_1'} \beta_{j_2}\beta_{j_2'}(\om_o)}}
 \beta_{j_2'}(\om_o)\\
&& \hspace*{-.8in}\times \EE\big[ {P}_{j_1l_1} \overline{{P}_{j_1'l_1'}} \overline{{P}_{j_2l_2}}{P}_{j_2'l_2'}(\om_o,\Lep;0)\big]
\exp\big( i (\beta_{j_1} -\beta_{j_1'}-\beta_{j_2}+\beta_{j_2'})(\omega_o+\omega) \Lep\big)  ,
\ea
which depends on the fourth-order moment of the propagator matrix
${\bf P}(\om,\Lep;0)$.
When scattering is strong, in the sense that $L$ is larger than the equipartition distance, these moments are \cite[Chapter 20]{fouque07}
\begin{eqnarray*}
&& \hspace*{-1.in}
\EE\big[ {P}_{j_1l_1} \overline{{P}_{j_1'l_1'}} \overline{{P}_{j_2l_2}}{P}_{j_2'l_2'}(\om_o,\Lep;0)\big]
 \approx
\left\{
\begin{array}{ll}
\frac{2}{N(N+1)}
&\mbox{ if } (j_1,l_1)=(j_1',l_1')=(j_2,l_2)=(j_2',l_2') \, ,\\
\frac{1}{N(N+1)}
&\mbox{ if }  (j_1,l_1)=(j_1',l_1') \neq (j_2,l_2)=(j_2',l_2') \, ,\\
\frac{1}{N(N+1)}
&\mbox{ if }  (j_1,l_1)=(j_2,l_2) \neq (j_1',l_1')=(j_2',l_2') \, ,\\
0 
& \mbox{ otherwise}\, .
\end{array}
\right.
\end{eqnarray*}
Therefore, we obtain that 
\ba
\nonumber
&& \hspace*{-1.in}
\EE\big[\cG_{11}(\omega,y,x)\cG_{12}(\omega,y',x)\big]
 =  \frac{1}{16 N(N+1)} 
\Big[\sum_{l=1}^N S_{ll} \beta_l(\om_o)  \Big]^2\\
&& \hspace*{1.in}
\times
\Big[
\Phi_{-1}(y,x)\Phi_0(y',x)+ \Phi_{-1}(y,y')\Phi_0(x,x)
\Big] ,
\ea
which gives the desired result.

\end{document}

%% file: mathmacros.tex
\DeclareMathAlphabet{\itbf}{OML}{cmm}{b}{it}
 \DeclareMathAlphabet\mathbfcal{OMS}{cmsy}{b}{n}

\renewcommand{\hat}{\widehat}
\newcommand\lb{\left<}
\newcommand\rb{\right>}
\def\EE{\mathbb{E}}
\def\RR{\mathbb{R}}

\def\bx{{{\itbf x}}}

\newcommand{\la}{\lambda}
\newcommand{\ep}{\varepsilon}

\newcommand{\om}{\omega}
\newcommand{\bka}{\boldsymbol{\kappa}}
\newcommand{\cI}{\mathscr{I}}
\newcommand{\cC}{\mathscr{C}}
\newcommand{\cG}{\mathscr{G}}
\newcommand{\cD}{\mathfrak{X}}
\newcommand{\cR}{\mathcal{R}}
\newcommand{\cL}{\mathcal{L}}
\newcommand{\cAd}{\mathcal{A}_{\rm d}}
\newcommand{\cAs}{\mathcal{A}_{\rm s}}
\newcommand{\Ir}{I^{({\rm r})}\hspace{-0.01in}}
\newcommand{\pr}{p^{({\rm r})}\hspace{-0.01in}}
\newcommand{\Iin}{I^{({\rm i})}}
\newcommand{\cIin}{\cI^{({\rm i})}}
\newcommand{\pin}{p^{({\rm i})}}
\newcommand{\hgr}{\hat{g}^{\,({\rm r})}\hspace{-0.01in}}
\newcommand{\hpin}{\hat{p}^{\,({\rm i})}}
\newcommand{\psc}{p^{({\rm s})}}
\newcommand{\hpsc}{\hat{p}^{\,({\rm s})}}
\newcommand{\cT}{\mathfrak{T}}

\newcommand{\Lep}{L_\ep}
\newcommand{\cLep}{\cL_\ep}

\newcommand{\bS}{{\bf S}}
\newcommand{\bD}{{\bf D}}
\renewcommand{\phi}{\varphi}